\documentclass{m2an}
\usepackage{amsmath,amsfonts,amssymb,latexsym,epic,eepic}
\usepackage[dvips]{graphics,epsfig,color}
\usepackage{psfrag,pstricks}

\newcommand{\dx}[1]{\, {\rm d}x}
\newcommand{\dedge}{\, {\rm d} \gamma}
\newcommand{\dt}{\delta t}
\newcommand{\dive}{\nabla\cdot}

\newcommand{\vs}{\mathrm{v}_{\edge,K}}
\newcommand{\vsp}{\mathrm{v}_{\edge,K}^+}
\newcommand{\vsm}{\mathrm{v}_{\edge,K}^-}

\newcommand{\edge}{\sigma}
\newcommand{\edged}{\varepsilon}
\newcommand{\edges}{{\cal E}}

\newcommand{\edgesint}{{\cal E}_{{\rm int}}}
\newcommand{\edgesext}{{\cal E}_{{\rm ext}}}
\newcommand{\mesh}{{\cal M}}

\newcommand{\fluxK}{F_{\edge,K}}
\newcommand{\fluxL}{F_{\edge,L}}
\newcommand{\fluxd}{F_{\edged,\edge}}
\newcommand{\normLdiscd}[2]{\hspace{.2em}|\hspace{-.1em}|#2|\hspace{-.1em}|_{h,#1}^2\hspace{.2em}}
\newcommand{\snormundisc}[2]{|#2|_{h,#1}}
\newcommand{\snormundiscd}[2]{|#2|_{h,#1}^2}
\newcommand{\norminf}[1]{\hspace{.2em}|\hspace{-.1em}| #1 |\hspace{-.1em}|_{\infty}\hspace{.2em}}
\newcommand{\normLd}[1]{\hspace{.2em}|\hspace{-.1em}| #1 |\hspace{-.1em}|_{\xLtwo(\Omega)}\hspace{.2em}}

\newcommand{\normHb}[1]{\hspace{.2em}|\hspace{-.1em}| #1 |\hspace{-.1em}|_{1,b}\hspace{.2em}}
\newcommand{\normHbd}[1]{\hspace{.2em}|\hspace{-.1em}| #1 |\hspace{-.1em}|_{1,b}^2\hspace{.2em}}

\newcommand{\norma}[1]{\hspace{.2em}|\hspace{-.1em}| #1 |\hspace{-.1em}|_{\rm a}\hspace{.2em}}
\newcommand{\normad}[1]{\hspace{.2em}|\hspace{-.1em}| #1 |\hspace{-.1em}|^2_{\rm a}\hspace{.2em}}
\newcommand{\normma}[1]{\hspace{.2em}|\hspace{-.1em}| #1 |\hspace{-.1em}|_{\rm -a}\hspace{.2em}}
\newcommand{\normmad}[1]{\hspace{.2em}|\hspace{-.1em}| #1 |\hspace{-.1em}|^2_{\rm -a}\hspace{.2em}}
\newcommand{\matB}{{\rm B}}
\newcommand{\matL}{{\rm L}}
\newcommand{\matM}{{\rm M}}
\newcommand{\matR}{{\rm R}}
\newcommand{\matQ}{{\rm Q}}

\newcommand{\grad}{\nabla}              

\def\eg{\emph{e.g.\/}}
\def\ie{\emph{i.e.\/}}

\def\etal{\emph{et al.\/}}

  \def\xR{\mathbb{R}}

\def\xLtwo{{\rm L}^{2}} 
\def\xLinfty{{\rm L}^{\infty}}

\def\xHone{{\rm H}^{1}}

\newcommand{\D}{\displaystyle}

\begin{document}

\title{An entropy preserving finite-element/finite-volume pressure correction scheme for the drift-flux model}

\author{L. Gastaldo}
\address{Institut de Radioprotection et de S\^{u}ret\'{e} Nucl\'{e}aire (IRSN) (laura.gastaldo@irsn.fr)}

\author{R. Herbin}
\address{Universit\'e de Provence, France (herbin@cmi.univ-mrs.fr)}

\author{J.-C. Latch\'e}
\address{Institut de Radioprotection et de S\^{u}ret\'{e} Nucl\'{e}aire (IRSN) (jean-claude.latche@irsn.fr)}

\begin{abstract}
We present in this paper a pressure correction scheme for the drift-flux model combining finite element and finite volume discretizations, which is shown to enjoy essential stability features of the continuous problem: the scheme is conservative, the unknowns are kept within their physical bounds and, in the homogeneous case (\ie\ when the drift velocity vanishes), the discrete entropy of the system decreases; in addition, when using for the drift velocity a closure law which takes the form of a Darcy-like relation, the drift term becomes dissipative.
Finally, the present algorithm preserves a constant pressure and a constant velocity through moving interfaces between phases.
To ensure the stability as well as to obtain this latter property, a key ingredient is to couple the mass balance and the transport equation for the dispersed phase in an original pressure correction step.
The existence of a solution to each step of the algorithm is proven; in particular, the existence of a solution to the pressure correction step is derived as a consequence of a more general existence result for discrete problems associated to the drift-flux model.
Numerical tests show a near-first-order convergence rate for the scheme, both in time and space, and confirm its stability.
\end{abstract}

%
%
%
%
%
%


\maketitle


\section{Introduction}

Dispersed two-phase flows and, in particular, bubbly flows are widely encountered in industrial applications as, for instance, nuclear safety studies, which are the context of the present work.
Within the rather large panel of models dealing with such flows, the simplest is the so-called drift-flux model, which consists in balance equations for an equivalent continuum representing both the gaseous and the liquid phase.
For isothermal flows, this approach leads to a system of three balance equations, namely the overall mass, the gas mass and the momentum balance, which reads:
\begin{equation}
\left|
\begin{array}{l} \displaystyle
\frac{\partial\,\rho}{\partial t}+\nabla\cdot(\rho\, u)=0
\\[2ex] \displaystyle
\frac{\partial\,\rho\,y}{\partial t}+\nabla\cdot(\rho\,y\, u)
    =-\nabla\cdot(\rho\,y\, (1-y)\, u_r) +\nabla \cdot (D \nabla y) 
\\[2ex] \displaystyle
\frac{\partial \rho\, u}{\partial t}+\nabla\cdot(\rho\, u\otimes u)+\nabla p-\nabla\cdot\tau(u)=f_v
\end{array}\right.
\label{cont_problem}\end{equation}
where $t$ stands for the time, $\rho$, $u$ and $p$ are the (average) density, velocity and pressure in the flow and $y$ stands for the gas mass fraction.
The diffusion coefficient $D$ represents in most applications small scale perturbations of the flow due to the presence of the dispersed phase, sometimes called "diphasic turbulence" and $u_r$ is the relative velocity between the liquid and the gaseous phase (the so-called drift velocity); for both these quantities, a phenomenologic relation must be supplied.
The forcing term $f_v$ may represent, for instance, the gravity forces.
The tensor $\tau$ is the viscous part of the stress tensor, given by the following expression:
\begin{equation}
\tau(u)=\mu \,(\grad u +\grad^t u)  - \frac{2}{3} \, \mu \, (\dive u)\ I
\label{tau}\end{equation}
For a constant viscosity, this relation yields:
\begin{equation}
\dive \tau = \mu \left[ \Delta u + \frac 1 3 \ \grad \dive u \right]
\label{tau2}\end{equation}
and, in this case, this term is dissipative (\ie\ for any regular velocity field $u$ vanishing on the boundary, the integral of $\dive \tau(u) \cdot u$ over the computational domain is non-negative).

\medskip
This system must be complemented by an equation of state, which takes the general form:
\begin{equation}
\rho=\varrho^{\,p,\alpha}(p,\alpha_g)=(1-\alpha_g) \rho_\ell+\alpha_g \varrho_g(p)
\label{state_law}
\end{equation}
where $\alpha_g$ stands for the void fraction and $\varrho_g(p)$ expresses the gas density as a function of the pressure; in the ideal gas approximation and for an isothermal flow, $\varrho_g(\cdot)$ is simply a linear function:
\begin{equation}
\varrho_g(p)=\frac{p}{a^2}
\label{gas_eos}\end{equation}
where $a$ is a constant characteristic of the gas, equal to the sound velocity in an isothermal (monophasic) flow.
The density of the liquid phase $\rho_\ell$ is assumed to be constant.
Introducing the mass gas fraction $y$ in \eqref{state_law} by using the relation $\alpha_g \, \varrho_g=\rho \, y$ leads to the following equation of state:
\begin{equation}
\rho=\varrho^{\,p,y}(p,\, y)= \frac{\varrho_g(p)\,\rho_\ell}{\rho_\ell\,y+(1-y)\ \varrho_g(p)}
\label{eos}\end{equation}

\medskip
The problem is supposed to be posed over $\Omega$, an open bounded connected subset of $\xR^d,\,d\leq 3$, and over a finite time interval $(0,T)$.
It must be supplemented by suitable boundary conditions, and initial conditions for $\rho$, $u$ and $y$.

\medskip
To design a numerical scheme for the solution of the system \eqref{cont_problem}, one is faced with several difficulties.
First, since the fluid density $\rho_\ell$ is supposed not to depend on the pressure, almost incompressible zones, \ie\ zones where the void fraction is low, may coexist in the flow with compressible zones, \ie\ zones where the void fraction remains significant.
This feature makes the problem particularly difficult to solve from a numerical point of view, because the employed numerical scheme will have to cope with a wide range of Mach numbers, starting from zero to, let us say, for low to moderate speed flows, a fraction of unity.
Second, the gas mass fraction $y$ can be expected, both for physical and mathematical reasons, to remain in the $[0,1]$ interval, and it appears strongly desirable that the numerical scheme reproduces this behaviour at the discrete level.
Finally, it appears from numerical experiments that, in order to avoid numerical instabilities, the algorithm should preserve a constant pressure through moving interfaces between phases (\ie\ contact discontinuities of the underlying hyperbolic system).
To obtain a scheme stable in the low Mach number limit, the solution that we adopt here is to use an algorithm inspired from the incompressible flow numerics, namely from the class of finite element pressure correction methods, and which degenerates to a classical projection scheme when the fluid density is constant.
The last two requirements are met thanks to an original pressure correction step in which the mass balance equation is solved simultaneously with a part of the gas mass balance.
For technical reasons, the solution of this latter equation is itself split in two steps, the first step thus being incorporated to the pressure correction step and the second one being performed independently.

\medskip
This work takes benefit of ideas developped in a wide literature, so we are only able to quote here some references, the choice of which will unfortunately probably appear somewhat arbitrary.
For a description of projection schemes for incompressible flow, see \eg\ \cite{gue-06-ove, mar-98-nav} and references herein.
An extension to barotropic Navier-Stokes equations close to the scheme developped here can be found in \cite{gal-07-an}, together with references to (a large number of) related works (see \eg\ \cite{har-71-num} for the seminal work and \cite{wes-01-pri} for a comprehensive introduction).
Extensions of pressure correction algorithms for multi-phase flows are more scarce, and seem to be restricted to iterative algorithms, often similar in spirit to the usual SIMPLE algorithm for incompressible flows \cite{spa-80-num,mou-03-pre,kuz-04-num}.
The gas mass balance equation, \ie\ the second equation of \eqref{cont_problem}, is a convection-diffusion equation which differs from the usual mass balance for chemical species in compressible multi-component flows studied by Larrouturou \cite{lar-91-how} by the addition of a non-linear term of the form $\dive \rho\, \varphi(y) \, u_r$, where $\varphi(\cdot)$ is a regular function such that $\varphi(0)=\varphi(1)=0$ (in the present case, $\varphi(y) = y \, (1-y)$).
In \cite{gas-07-on}, we propose a finite-volume scheme for the numerical approximation of this type of equation, and we prove the existence and uniqueness of the solution, together with the fact that it remains within physical bounds, \ie\ within the interval $[0,1]$.
Here, the proof of the same results combines arguments from both \cite{lar-91-how} and \cite{gas-07-on}.

\medskip
Several theoretical issues concerning the proposed scheme are studied in this paper.
First, the existence of a solution to the pressure correction step, which consists in an algebraic non-linear system, is obtained by a topological degree argument.
Second, we address the stability of the scheme.
At the continuous level, the existence of an entropy for the system when the drift velocity vanishes (\ie\ the homogeneous model) is well-known.
In addition, it is shown in \cite{gui-07-ada}, by a Chapman-Enskog expansion technique, that the two-fluid model can be reduced to the drift-flux model when a strong coupling of both phases is assumed, with a Darcy-like closure relation for the drift velocity, \ie\ an expression of the form:
\begin{equation}
u_r =  \frac 1 \lambda\ (1-\alpha_g)\,\alpha_g\,\frac{\varrho_g(p)-\rho_\ell}{\rho}\ \grad p
\label{u_r} \end{equation}
where $\lambda$ is a positive phenomenological coefficient.
The same relation can also be obtained by neglecting in the two-fluid model the difference of acceleration between both phases \cite{sok-04-sim}.
With such an expression for $u_r$, the drift term becomes a second order term, and it is shown in \cite{gui-07-ada} that it is consistent with the entropy of the homogeneous model (\ie\ that it generates a non-negative dissipation of the entropy).
These results are proven here at the discrete level: up to a minor modification of the proposed scheme, which seems useless in practice, the entropy is conserved when $u_r$ is equal to zero, and when the closure relation \eqref{u_r} applies and with a specific discretization, the drift term generates a dissipation.

\medskip
This paper is built as follows.
The fractional step algorithm for the solution of the whole problem is first presented in section \ref{sec:algo_z}, together with some of its properties: the existence of a solution to each step of the algorithm, the fact that the unknowns are kept within their physical bounds and that the algorithm is able to preserve a constant pressure and a constant velocity through moving interfaces between phases.
The proof of the existence of the solution to the pressure correction step is obtained as a consequence of a more general existence theory for some discrete problems associated to the drift-flux model, which is exposed in the appendix.
Next two sections are devoted to the stability analysis of the scheme; after establishing estimates for the work of the pressure forces (section \ref{sec:p_work}), we first address the case $u_r=0$ (section \ref{sec:u_r=0}), then the case where $u_r$ is given by the Darcy-like closure relation \eqref{u_r} (section \ref{sec:u_r=darcy}).
Finally, numerical tests are reported in section \ref{sec:numerical_z}; they include a problem exhibiting an analytical solution which allows to assess convergence properties of the discretization, a sloshing transient in a cavity, and the evolution of a bubble column.

\medskip
For the sake of simplicity, we suppose for the presentation of the scheme and its analysis (sections \ref{sec:algo_z}, \ref{sec:p_work} and \ref{sec:stab_z}) that the velocity is prescribed to zero on the whole boundary $\partial \Omega$ of the computational domain, and that the gas mass flux through $\partial \Omega$, so both $u_r$ and the component of $\grad y$ normal to the boundary, also vanishes.
Moreover, the analysis of the scheme assumes that pure liquid zones does not exist in the flow, which, with the proposed algorithm, is a consequence that such zones are not present at the initial time (\ie, at $t=0$, $y \in (0,1]$); getting rid of this latter limitation at the theoretical level seems indeed to be a difficult task.
However, the numerical tests presented in section \ref{sec:numerical_z} are not restricted to theses situations.
In particular, $y=0$ in the liquid column in the sloshing problem, up to spurious phases mixing by the numerical diffusion near the free surface; it is also the case at the initial time in the bubble column simulation.

\medskip
In the presentation of the scheme, the drift velocity is supposed to be known, \ie\ to be given by a closure relation independent of the unknowns of the problem, and this still holds in numerical experiments.
The case where $u_r$ is given by \eqref{u_r} is thus only treated from a theoretical point of view in section \ref{sec:u_r=darcy}.


\section{The numerical algorithm}\label{sec:algo_z}

We present in this section the numerical scheme considered in this paper.
We begin by describing the proposed scheme in the time semi-discrete setting, then we introduce the spatial discretization spaces and we detail the discrete approximation and the properties for each step of the algorithm at hand.


\subsection{Time semi-discrete formulation}

Let us consider a partition $0=t_0 < t_1 <\ldots < t_N=T$ of the time interval $(0,T)$, which is supposed uniform for the sake of simplicity.
Let $\dt$ be the constant time step $\dt=t_{n+1}-t_n$ for $n=0,1,\ldots,N$. 
In a time semi-discrete setting, the algorithm proposed in this paper is the following three steps scheme:
\begin{itemize}
\item[1 -] solve for $\tilde u^{n+1}$
\begin{equation}
\frac{\rho^n\ \tilde u^{n+1}-\rho^{n-1}\ u^n}{\dt}+\dive(\rho^n \ u^n \otimes\tilde u^{n+1})+\nabla p^n -\dive \tau(\tilde u^{n+1})= f_v^{n+1}
\label{qdm}
\end{equation}
\item[2 -] solve for $p^{n+1}$, $u^{n+1}$, $\rho^{n+1}$ and $z^{n+1}$
\begin{equation}
\left| \begin{array}{l} \displaystyle
\rho^n \ \frac{u^{n+1}-\tilde u^{n+1}}{\dt}+\nabla(p^{n+1}-p^n)=0
\\[2ex] \displaystyle
\frac{\varrho^{\,p,z}(p^{n+1},\ z^{n+1})-\rho^n}{\dt}+\dive(\varrho(p^{n+1},\ z^{n+1})\ u^{n+1})=0
\\[2ex] \displaystyle
\frac{z^{n+1}-\rho^n y^n}{\dt}+\dive(z^{n+1} \ u^{n+1})=0
\\[2ex] \displaystyle
\rho^{n+1}=\varrho^{\,p,z}(p^{n+1},z^{n+1})
\end{array}\right.
\label{projection}\end{equation}
\item[3 -] solve for $y^{n+1}$
\begin{equation}
\frac{\rho^{n+1}y^{n+1}-z^{n+1}}{\dt}+\dive(\rho^{n+1}\ y^{n+1}(1-y^{n+1})\ u_r^{n+1})=\dive( D\grad y^{n+1})
\label{splitting}
\end{equation}
\end{itemize}

The first step consists in a classical semi-implicit solution of the momentum balance equation to obtain a predicted velocity.

\medskip
Step 2 is an original nonlinear pressure correction step, which couples the mass balance equation (second equation) with the transport terms of the gas mass balance equation (third equation).
A new unknown is introduced in this step instead of the gas mass fraction, the partial gas density $z$ given by $z=\rho\,y$.
Thus, the equation of state must be reformulated to express the mixture density as a function of the partial gas density and of the pressure, which, from equation \eqref{eos}, yields:
\begin{equation}
\rho=\varrho^{\,p,z}(p,z)=z\left(1-\frac{\rho_\ell\,a^2}{p}\right)+\rho_\ell
\label{eos_z}
\end{equation}
When the liquid and the gas densities are very different, this law presents much less steep variations than the relation linking the density and the mass fraction $y$, specially in the neighbourhood of $y=0$; this change of variable thus makes the resolution of this step much easier, and the overall algorithm more robust.
In counterpart, it leads to split the gas mass balance equation: transport terms are dealt with in the present step, and the gas mass fraction is corrected in a next step (step 3) to take into account the drift terms.
The pressure correction step would degenerate in the usual projection step as used in incompressible flows solvers if the density was constant (\ie\ $z=0$).
Taking (at the algebraic level, see section \ref{sec:p_correction}) the divergence of the first relation of \eqref{projection} and using the second one to eliminate the unknown velocity $u^{n+1}$ yields a non-linear elliptic problem for the pressure. 
Solving at the same time this elliptic problem and the third equation by a Newton's algorithm, we obtain the pressure and the gas mass fraction.
Once the pressure is computed, the first relation yields the updated velocity and the fourth one gives the end-of-step density.

\medskip
Finally, in the third step, the remaining terms of the gas mass balance are considered, and the end-of-step gas mass fraction is computed.

\medskip
The motivations of this time discretization are the following ones: to keep the mass fraction $y$ in the physical range $[0,1]$, to allow the transport of phases interfaces without generating spurious pressure and velocity variations and to ensure the stability of (\ie\ the conservation of the entropy by) the scheme.
To show how this time splitting algorithm achieves these goals is the aim of the remaining of this paper.


\subsection{Spatial discretization}\label{sec:spacedisc}

Let $\mesh$ be a decomposition of the domain $\Omega$ either into convex quadrilaterals ($d=2$) or hexahedrons ($d=3$) or in simplices. 
By $\edges$ and $\edges(K)$ we denote the set of all $(d-1)$-edges $\edge$ of the mesh and of the element $K \in \mesh$ respectively.
The set of edges included in the boundary of $\Omega$ is denoted by $\edgesext$ and the set of internal ones (\ie\ $\edges \setminus \edgesext$) is denoted by $\edgesint$.
The decomposition $\mesh$ is supposed to be regular in the usual sense of the finite element literature (e.g. \cite{cia-91-bas}), and, in particular, $\mesh$ satisfies the following properties: $ \bar\Omega=\bigcup_{K\in \mesh} \bar K$; if $K,\,L \in \mesh,$ then $\bar K \cap \bar L=\emptyset$ or $\bar K\cap \bar L$ is a common edge of $K$ and $L$, which is denoted by $K|L$.
For each internal edge of the mesh $\edge=K|L$, $n_{KL}$ stands for the normal vector of $\edge$, oriented from $K$ to $L$.
By $|K|$ and $|\edge|$ we denote the measure, respectively, of $K$ and of the edge $\edge$.

\medskip
For stability reasons, the spatial discretization must preferably be based on pairs of velocity and pressure approximation spaces satisfying the so-called {\it inf-sup} or Babuska-Brezzi condition (\eg\ \cite{bre-91-mix}).
Among these elements, nonconforming approximations with degrees of freedom for the velocity located at the center of the faces seem to be well suited to a coupling with a finite volume treatment of the other equations, as is proposed hereafter for the gas mass balance; this is the choice made here.
The spatial discretization thus relies either on the so-called "rotated bilinear element"/$P_0$ introduced by Rannacher and Turek \cite{ran-92-sim} for quadrilateral or hexahedric meshes, or on the Crouzeix-Raviart element (see \cite{cro-73-con} for the seminal paper and, for instance, \cite[p. 83--85]{ern-05-aid} for a synthetic presentation) for simplicial meshes.
The reference element $\widehat K$ for the rotated bilinear element is the unit $d$-cube (with edges parallel to the coordinate axes); the discrete functional space on $\widehat K$ is $\tilde{Q}_{1}(\widehat K)^d$, where $\tilde{Q}_{1}(\widehat K)$ is defined as follows:
\[
\tilde{Q}_{1}(\widehat K)= {\rm span}\left\{1,\,(x_{i})_{i=1,\ldots,d},\,(x_{i}^{2}-x_{i+1}^{2})_{i=1,\ldots,d-1}\right\}
\]
The reference element for the Crouzeix-Raviart is the unit $d$-simplex and the discrete functional space is the space $P_1$ of affine polynomials.
For both velocity elements used here, the degrees of freedom are determined by the following set of nodal functionals:
\begin{equation}
\D \left\{F_{\edge,i},\ \edge \in \edges(K),\, i=1,\ldots,d\right\}, 
\qquad F_{\edge,i}(v)=|\edge|^{-1}\int_{\edge} v_{i} \dedge
\label{vdof}\end{equation}
The mapping from the reference element to the actual one is, for the Rannacher-Turek element, the standard $Q_1$ mapping and, for the Crouzeix-Raviart element, the standard affine mapping.
Finally, in both cases, the continuity of the average value of discrete velocities (\ie, for a discrete velocity field $v$, $F_{\sigma,i}(v),\ 1\leq i \leq d$) across each face of the mesh is required, thus the discrete space $W_{h}$ is defined as follows:
\[
\begin{array}{ll} \displaystyle
W_h =
& \displaystyle
\lbrace
\ v_h\in L^{2}(\Omega)\,:\, v_h|_K \in\tilde{Q}_{1}(K)^d,\,\forall K\in \mesh;
\\[1ex] & \D \hspace*{2cm}
\ F_{\edge,i}(v_h) \mbox{ continuous across each edge } \sigma \in {\edgesint},\mbox{ for } 1\leq i \leq d \,;
\\[1ex] & \D \hspace*{2cm}
\ F_{\edge,i}(v_h)=0,\ \forall \edge \in \edgesext,\ 1\leq i \leq d\ \rbrace
\end{array}
\]
For both Rannacher-Turek and Crouzeix-Raviart discretizations, the pressure is approximated by the space $L_{h}$ of piecewise constant functions:
\[
L_h=\left\{q_h\in L^{2}(\Omega)\,:\, q_h|_K=\mbox{ constant},\,\forall K\in \mesh\right\}
\]
Since only the continuity of the integral over each edge of the mesh is imposed, the velocities are discontinuous through each edge; the discretization is thus nonconforming in $H^1(\Omega)^d$.
These pairs of approximation spaces for the velocity and the pressure are \textit{inf-sup} stable, in the usual sense for "piecewise $\xHone$" discrete velocities, \ie\ there exists $c_{\rm i}>0$ independent of the mesh such that:
\[
\forall p \in L_h, \qquad \sup_{v \in W_h} \frac{\D \int_{\Omega,h} p \, \nabla \cdot v \, {\rm d}x}{\normHb{v}} \geq c_{\rm i} \normLd{p-m(p)}
\]
where $m(p)$ is the mean value of $p$ over $\Omega$, the symbol $\D \int_{\Omega,h}$ stands for $\D \sum_{K\in\mesh} \int_K$ and $\normHb{\cdot}$ stands for the broken Sobolev $\xHone$ semi-norm:
\[
\normHbd{v}=\sum_{K\in \mesh} \int_K |\nabla v |^2 \, {\rm d}x=\int_{\Omega,h}| \nabla v |^2 \, {\rm d}x
\]

\medskip
From the definition \eqref{vdof}, each velocity degree of freedom can be univoquely associated to an element edge.
Hence, the velocity degrees of freedom may be indexed by the number of the component and the associated edge, and the set of velocity degrees of freedom reads:
\[
\lbrace v_{\edge,i},\ \edge \in \edgesint,\ 1 \leq i \leq d \rbrace
\]
We define $v_\edge=\sum_{i=1}^d v_{\edge,i}\, e^{(i)}$ where $e^{(i)}$ is the $i^{th}$ vector of the canonical basis of $\xR^d$.
We denote by $\varphi_{\edge}^{(i)}$ the vector shape function associated to $v_{\edge,i}$, which, by the definition of the considered finite elements, reads:
\[
\varphi_{\edge}^{(i)}=\varphi_\edge \, e^{(i)}
\]
where $\varphi_\edge$ is a scalar function.

\medskip
Each degree of freedom for the pressure is associated to a mesh $K$, and the set of pressure degrees of freedom is denoted by $\lbrace p_K,\ K \in \mesh \rbrace$.
As the pressure, the density $\rho$, the gas mass fraction $y$ and the gas partial density $z$ are approximated by piecewise constant functions over each element, and the associated sets of degrees of freedom are denoted by $\lbrace \rho_K,\ K \in \mesh \rbrace$, $\lbrace y_K,\ K \in \mesh \rbrace$ and $\lbrace z_K,\ K \in \mesh \rbrace$ respectively.


\subsection{Spatial discretization of the momentum balance equation}\label{sec:momentum}

The main difficulty in the discretization of the momentum balance equation is to build a discrete convection operator which enjoy the discrete analogue of the kinetic energy relation, that is:
\[
\int_\Omega \left[ \frac{\partial \rho\, u}{\partial t}+\nabla\cdot(\rho\, u\otimes u) \right] \cdot u \ {\rm d}x=
\frac d {dt} \int_\Omega \frac 1 2 \,\rho \,|u|^2 
\qquad \mbox{provided that} \quad \frac{\partial \rho}{\partial t}+\nabla\cdot(\rho\, u)=0
\]
To this purpose, we follow an idea developped in \cite{bab-07-an}, and already exploited for the same problem as here in \cite{gas-07-on}.
The idea is to derive a finite-volume-like discretization of the convection operator, in order to apply the following result \cite{gal-07-an}.

\begin{thrm}[Stability of a finite volume advection operator]\label{VF1_2}Let $(\rho_K^\ast)_{K\in \mesh}$ and $(\rho_K)_{K\in \mesh}$ be two families of positive real numbers satisfying the following set of equation:
\begin{equation}
\forall K \in \mesh,\qquad \frac{|K|}{\dt} \ (\rho_K - \rho^\ast_K) + \sum_{\edge=K|L} \fluxK=0
\label{mass_bal}\end{equation}
where $\fluxK$ is a quantity associated to the edge $\edge$ and to the control volume $K$; we suppose that, for any internal edge $\edge=K|L$, $\fluxK=-\fluxL$.
Let $(s_K^\ast)_{K\in \mesh}$ and $(s_K)_{K\in \mesh}$ be two families of real numbers.
The following stability property holds:
\begin{equation}
\hspace*{-2ex}\sum_{K \in \mesh} z_K \left[ \frac{|K|}\dt (\rho_K s_K -\rho_K^\ast s_K^\ast)+
\sum_{\edge=K|L} \fluxK\ \frac{s_K+s_L} 2 \right]
\geq
\frac{1}{2} \sum_{K\in\mesh} \frac{|K|}\dt \left[ \rho_K s_K^2 -\rho_K^\ast {s_K^\ast}^2\right]
\label{ecinetique} \end{equation}
\end{thrm}

\begin{figure}[htb]
\begin{tabular}{lr}
\begin{minipage}{0.35\textwidth}
\psfrag*{K}{$K$}
\psfrag*{D}{$D_{K,\edge}$}
\psfrag*{sigma}{$\edge\in\edges(K)$}
\psfrag*{sigmad}{$\edged\in\edges(D_\edge)$}
\scalebox{0.75}{\includegraphics*{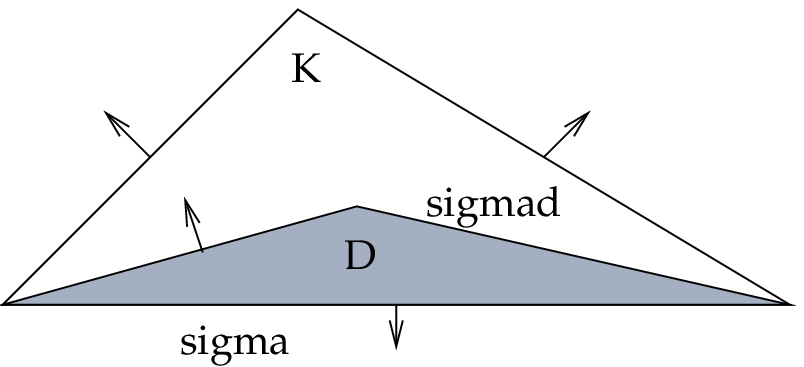}}
\end{minipage}
&
\begin{minipage}{0.40\textwidth}
\psfrag*{K}{$K$}
\psfrag*{D}{$D_{K,\edge}$}
\psfrag*{sigma}{$\edge\in\edges(K)$}
\psfrag*{sigmad}{$\edged\in\edges(D_\edge)$}
\scalebox{0.75}{\includegraphics*{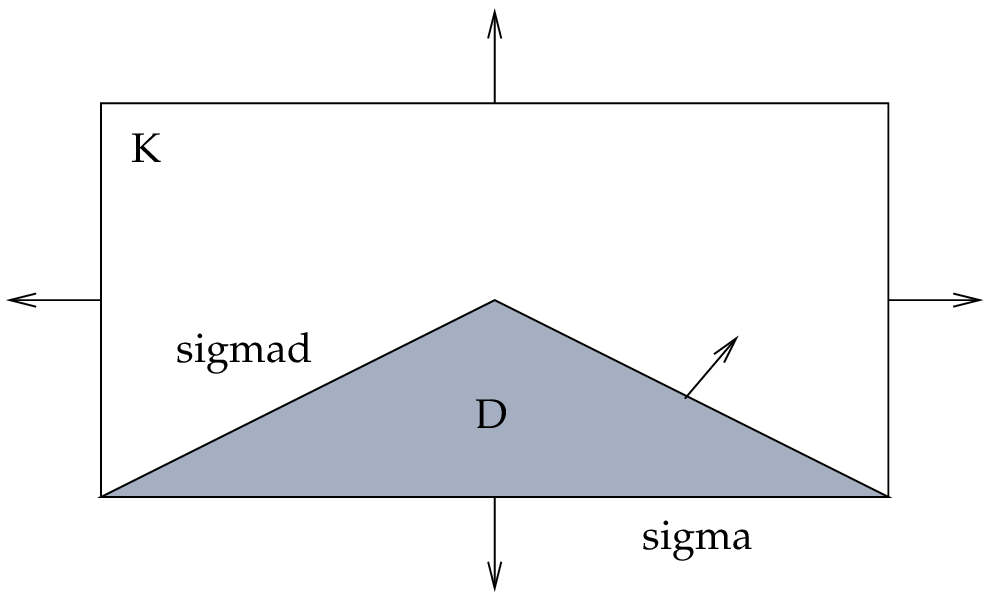}}
\end{minipage}
\end{tabular}
\caption{Diamond-cells for the Crouzeix-Raviart and Rannacher-Turek element. \label{diamonds}}
\end{figure}

\medskip
To this purpose, we first define a control volume for each degree of freedom of the velocity, that is, in view of the discretization used here, around each barycenter of an internal edge.
Let $\edge=K|L$ and $D_{K,\edge}$ be the conic volume having $\edge$ for basis and the mass center of $K$ as additional vertex (see figure \ref{diamonds}). 
The volume $D_\edge=D_{K,\edge}\cup D_{L,\edge}$ is referred to as the "diamond cell" associated to $\edge$ and $D_{K,\edge}$ is the half-diamond cell associated to $\edge$ and $K$.
For the Crouzeix-Raviart element and, for the Rannacher-Turek element, when the mesh is a rectangle (in two dimensions) or a cuboid (in three dimensions), the integral of the shape function associated to the edge $\edge$ over the element $K$ is the measure of the half-diamond cell $D_{K,\edge}$.
Thus, the application of the mass lumping to the terms of the form $\rho u$ leads, in the equations associated to the velocity on the edge $\edge$, to a discrete expression of the form $\rho_\edge u_\edge$, where $\rho_\edge$ results from an average of the values taken by the density in the two elements adjacent to $\edge$, weighted by the measure of the half-diamonds:
\begin{equation}
\forall \edge \in \edgesint,\qquad |D_\edge|\ \rho_\edge= |D_{K,\edge}|\ \rho_K + |D_{L,\edge}|\ \rho_L
\label{rho_faces2}\end{equation}
where $|D_\edge|$ is the measure of the diamond cell $D_\edge$, $|D_{K,\edge}|$ and $|D_{L,\edge}|$ are the measure of the half-diamond cells associated respectively to $\edge$ and $K$ and to  $\edge$ and $L$.
This lumped time derivative term naturally combines with a discretization of the advective term of the form:
\[
\forall \edge \in \edgesint, \qquad \mbox{ advection term } \sim
\sum_{\edged \in \edges(D_\edge)} \fluxd\ u_\edged
\]
where $\edges(D_\edge)$ is the set of the edges of $D_\edge$, $u_\edged$ is a centered approximation of $u$ on $\edged \in \edges(D_\edge)$ and $\fluxd$ is a mass flux through $\edged$.
To proceed, we must now derive for this latter quantity an approximation which satisfies the compatibility condition \eqref{mass_bal} of theorem \ref{VF1_2} (in fact, the discrete mass balance over the diamond cells).
Suppose that we are able to build, for any control volume $K$, a field $\widetilde{\rho u}_K(x)$ such that $\nabla \cdot \widetilde{\rho u}_K(x)$ remains constant inside the element $K$ and that we take for the mass flux $\fluxd$ through each diamond cell edge $\edged$ included in $K$:
\[
\fluxd=\int_\edged \widetilde{\rho u}_K(x) \cdot n_{\edged,\edge} \ {\rm d}\gamma (x)
\]
where $n_{\edged,\edge}$ is the normal vector to $\edged$ outward $D_\edge$. 
As the divergence of $\widetilde{\rho u}_K$ is constant over $K$, it may be checked that, if the flux of $\widetilde{\rho u}_K$ through each edge of $K$ is the same as the mass flux used in a discrete mass balance over K, let say $|\edge|\, (\rho \ u)_\edge\cdot n_\edge$, this mass balance "carry over" the half-diamond cells $D_{K,\edge}$, which, by summation over the two half-diamond cells, yields a compatibility condition of the desired form \cite{gas-07-on}.
Such a field $\widetilde{\rho u}_K(x)$ is derived for the Crouzeix-Raviart element by direct interpolation (\ie\ using the standard expansion of the Crouzeix-Raviart elements) of the quantities $((\rho \ u)_\edge)_{\edge\in \edges(K)}$:
\[
\widetilde{\rho u}_K(x)=\sum_{\edge\in\edges(K)}\varphi_\edge(x) \ (\rho \ u)_\edge
\]
For the Rannacher-Turek element, when the mesh is a rectangle or a cuboid, it is obtained by the following interpolation formula:
\[
\widetilde{\rho u}_K(x)=\sum_{\edge\in\edges(K)}\alpha_\edge(x\cdot n_\edge)\ \left[(\rho \ u)_\edge\cdot n_\edge\right]\ n_\edge
\]
where the $\alpha_\edge(\cdot)$ are affine interpolation functions which are determined in such a way that the desired conditions hold, \ie\ that the flux of $\widetilde{\rho u}$ through each edge $\edge$ of $K$ is $|\edge|\, (\rho \ u)_\edge\cdot n_\edge$.
Extension to more general grids is underway.
Finally, since, in the proposed fractional step algorithm, the mass balance equation is considered only when the solution of the momentum balance is achieved, to obtain the desired compatibility condition \eqref{mass_bal}, we use the mass balance at the previous time step: the approximations of the density in the time derivative term are shifted of one time step and the quantities $((\rho \ u)_\edge)_{\edge\in \edges(K)}$ used to compute the mass fluxes $\fluxd^n$ are chosen to be the mass fluxes obtained in the discrete mass balance at the previous time step.
Since standard finite elements techniques are used to discretize the term $\nabla p^n - \dive \tau(\tilde u^{n+1})$, this yields the following discrete momentum balance equation:
\begin{equation}
\begin{array}{l} \displaystyle
\forall \edge \in \edgesint,\mbox{ for } 1 \leq i \leq d,
\\[3ex] \displaystyle \hspace{5ex}
\frac{|D_\edge|} \dt\ (\rho_\edge^n\ \tilde u_{\edge,i}^{n+1}- \rho^{n-1}_\edge \ u_{\edge,i}^n)
+ \sum_{\stackrel{\scriptstyle \edged \in \edges(D_\edge),}{\scriptstyle \edged=D_\edge |  D_{\edge'}}}\frac 1 2\ \fluxd^n\ (\tilde u_{\edge,i}^{n+1} + \tilde u_{\edge',i}^{n+1})
\\ \hspace{40ex}\displaystyle
+ a_d(\tilde u^{n+1},\varphi_\edge^{(i)})
-\int_{\Omega,h}  p^n\ \nabla \cdot \varphi_\edge^{(i)}
= \int_\Omega f^{n+1} \cdot \varphi_\edge^{(i)}
\end{array}
\label{momentum_2}
\end{equation}
where, the bilinear form $a_d(\cdot,\cdot)$ represents the viscous term and, $\forall v \in W_h,\ \forall w \in W_h$, $a_d(v,w)$ is defined as follows:
\[
a_d(v,w)=
\left| \begin{array}{ll} \displaystyle
\mu \int_{\Omega,h} \left[ \nabla v : \nabla w + \frac 1 3 \, \dive v\ \dive w \right]\ {\rm d}x \qquad
&
\mbox{if \eqref{tau2} holds (case of constant viscosity),}
\\[3ex] \displaystyle
\int_{\Omega,h} \tau(v) : \nabla w\ {\rm d}x
&
\mbox{with $\tau$ given by \eqref{tau} otherwise.}
\end{array}\right.
\]

Note that, for Crouzeix-Raviart elements, a combined finite volume/finite element method similar to the technique employed here has already been analysed for a transient non-linear convection-diffusion equation by Feistauer and co-workers \cite{ang-98-ana, dol-02-err, fei-03-mat}.

\medskip 
As a consequence of the stability of the convection operator, we have the following regularity result.

\begin{lmm}[Properties of the numerical scheme - velocity prediction]
Let us assume that the viscous term is dissipative (\ie\ $\forall v \in W_h, \ a_d(v,v)\geq 0$, which holds for the form of $a_d(\cdot,\cdot)$ used in case of a constant viscosity); then the first step of the scheme, namely the velocity prediction step, has a unique solution.
\label{exist_step_v}\end{lmm}

\begin{rmrk}[First time step]
To ensure the compatibility condition \eqref{mass_bal} at first time step, a prediction step must be used to initialize the density:
\begin{equation}
\frac{\rho^0-\rho^{-1}}{\dt}+\dive(\rho^0 u^{-1})=0
\label{pred_rho}
\end{equation}
where $\rho^{-1}$ and $u^{-1}$ are suitable approximations for the initial density and the velocity, respectively.
\label{first_step}\end{rmrk}


\subsection{Spatial discretization of the pressure correction step}\label{sec:p_correction}

The discretization of the first equation of the pressure correction step is consistent with the momentum balance one, \ie\ we use a mass lumping technique for the unsteady term and a standard finite element formulation for the gradient of the pressure increment:
 \[
\forall \edge \in \edgesint,\mbox{ for } 1 \leq i \leq d,\qquad
\frac{|D_\edge|}{\dt} \, \rho_\edge^{n}\,(u_{\edge,i}^{n+1} - \tilde u_{\edge,i}^{n+1}) 
+ \int_{\Omega,h} (p^{n+1}-p^{n})\ \nabla \cdot \varphi_{\edge}^{(i)} \ dx = 0
\]
As the pressure is piecewise constant, the transposed of the discrete gradient operator takes the form of the finite volume standard discretization of the divergence based on the finite element mesh, thus the previous relation can been rewritten as follows:
\begin{equation}
\forall \edge\in\edgesint, \ \edge=K|L,\qquad
|\frac{|D_\edge|}{\dt} \, \rho_\edge^{n}\,(u_\edge^{n+1}-\tilde u_\edge^{n+1}) 
+ |\edge| \left[ (p_K^{n+1}-p_K^{n})-(p_L^{n+1}-p_L^{n}) \right] n_{KL} = 0
\label{pc1}\end{equation}

\medskip
Similarly, as the density is piecewise constant, the approximation of the time derivative of the density in the mass balance will also look as a finite volume term.
This point suggests a finite volume discretization of this latter equation, which reads:
\begin{equation}
\forall K \in \mesh, \qquad
\frac{|K|}{\dt} \ (\varrho^{\,p,z}(p^{n+1}_K, z^{n+1}_K)-\rho^n_K)
+ \sum_{\edge=K|L}\fluxK^{n+1}=0
\label{pc2}\end{equation}
To ensure the positivity of the density, we use an upwinding technique for the convection term, then the mass flux from $K$ across $\edge=K|L$, $\fluxK^{n+1}$, is expressed as follows:
\[
\fluxK^{n+1}=|\edge|\ (\rho^{n+1} u^{n+1})_{|\edge}\cdot n_{\edge}
=(\vsp)^{n+1}\, \varrho^{\,p,z}(p_K^{n+1}, z^{n+1}_K) - (\vsm)^{n+1}\, \varrho^{\,p,z}(p_L^{n+1}, z^{n+1}_L)
\]
where $(\vsp)^{n+1}$ and $(\vsm)^{n+1}$ stands respectively for $\max (\vs^{n+1},\ 0)$ and $-\min (\vs^{n+1},\ 0)$ with $\vs^{n+1}=|\edge|\, u_\edge^{n+1} \cdot n_{KL}$.

\medskip
Consistently with the mass balance equation, we use for the discretization of the third relation of \eqref{projection}, \ie\ the transport of the gas partial density $z$,  a finite volume method with an upwind technique for the convection term $\dive( z \,u)$.
This yields the following discrete equation:
\begin{equation}
\forall K \in \mesh, \qquad
\frac{|K|}{\dt}\ (z^{n+1}-\rho_{K}^{n}\,y_{K}^n)
+\sum_{\edge=K|L}(\vsp)^{n+1}\, z_K^{n+1}- (\vsm)^{n+1}\, z^{n+1}_L=0
\label{pc3}\end{equation}

In the following lemma, we state some properties of this pressure correction step which are obtained as a particular case of the existence theory presented in section \ref{sec:existence_z}.

\begin{lmm}[Properties of the numerical scheme - pressure correction step]
Let the density of the liquid phase be constant and the gas phase obeys the ideal gas law.
Then, under the assumption that, $\forall K \in \mesh,\ \rho^n_K>0$ and $y^n_K \in (0,1]$, the system \eqref{pc1}-\eqref{pc3} has a solution, and any solution of this step is such that:
\[
\forall K\in\mesh, \qquad \rho_K^{n+1} > 0, \quad  p_K^{n+1} >0, \quad z_K^{n+1}>0 \quad \mbox{and} \quad
\frac{z_K^{n+1}}{\rho_K^{n+1}} \in (0,1]
\] 
\label{exist_step_p}\end{lmm}

Let us now turn to the practical solution of this pressure correction step.
Keeping the same notation for the unknown functions and the vectors gathering their degrees of freedom, the algebraic formulation of this step reads:
\begin{equation}
\begin{array}{l} 
\left| \begin{array}{l} \displaystyle
\frac{1}{\dt} \matM_{\rho^{n}} \, (u^{n+1} - \tilde u^{n+1})+ \matB^t\, (p^{n+1} -p^{n})=0
\\[2ex] \displaystyle
\frac{1}{\dt} \matR\, (\varrho^{\,p,z}(p^{n+1},\,z^{n+1})-\rho^n) - \matB \matQ_{\rho^{n+1}}^{\rm up} u^{n+1}=0
\\[2ex] \displaystyle
\frac{1}{\dt} \matR\, (z^{n+1}-\rho^n\,y^n)- \matB \matQ_{z^{n+1}}^{\rm up} u^{n+1}=0
\end{array} \right .
\end{array}
\label{proj-alg3}\end{equation}
In the first relation, ${\rm M}_{\rho^{n}}$ stands for the diagonal mass matrix weighted by the density at $t^{n}$ (at edge center) $\rho^n_\edge$, so the diagonal entry of ${\rm M}_{\rho^{n}}$ associated to the internal edge $\edge$ and the component $i$ reads $(\matM_{\rho^{n}})_{\edge,i}= |D_\edge|\,\rho_\edge^{n}$.
The matrix $\matB^t$ of $\xR^{N \times M}$, where $N=d\ {\rm card}\,(\edgesint)$ and $M={\rm card}\,(\mesh)$, is associated to the gradient operator; consequently, the matrix $\matB$ is associated to the opposite of the divergence operator.
In the second and in the third relation, ${\rm Q}_{w^{n+1}}^{\rm up}$ (with $w=\rho$ or $w=z$) is a diagonal matrice, the entry of which corresponding to an edge $\edge\in \edgesint,\ \edge=K|L$, is obtained by just taking $w$ at $t^{n+1}$ in the element located upstream of $\edge$ with respect to $u^{n+1}$, {\it i.e.} either $w_K^{n+1}$ or $w_L^{n+1}$.
The matrix ${\matR}$ is diagonal and, for any $K\in\mesh$, its entry $\matR_K$ is the measure of the element $K$.

\medskip
The elliptic problem for the pressure is obtained by multiplying the first relation of (\ref{proj-alg3}) by $\matB\ \matQ_{\rho^{n+1}}^{\rm up}\ (\matM_{\tilde \rho^{n+1}})^{-1}$ and using the second one.
This equation reads:
\begin{equation}
\matL\, p^{n+1}+ \frac{1}{\dt^{2}}\matR\, \varrho^{\,p,z}(p^{n+1},\,y^{n+1}) =
\matL\,p^{n}+\frac{1}{\dt^{2}}\matR\, \rho^n +\frac{1}{\dt}\matB\,\matQ_{\rho^{n+1}}^{\rm up}\,\tilde u^{n+1}
\label{step31_3}\end{equation}
where, as seen in \cite{gal-07-an}, $\matL=\matB\,\matQ_{\rho^{n+1}}^{\rm up}\,(\matM_{\rho^{n}})^{-1}\,\matB^t$ can be equivalently evaluated in the "finite volume way" by the following relation, valid for each element $K$:
\[
(\matL\, p^{n+1})_K=\sum_{\edge=K|L} \ 
\frac{\D \rho_{{\rm up},\edge}^{n+1}}{\rho^{n}_\edge}\ \frac{|\edge|^2}{|D_\edge|}\, (p_K-p_L)
\]
where $\rho_{{\rm up},\edge}$ stands for the upwind density associated to the edge $\edge$.
One recognize in this relation a usual finite volume diffusion operator, with a particular diffusion coefficient which, for instance, can be evaluated for rectangular parallelepipedic control volumes as $d \ \rho_{{\rm up},\edge}^{n+1}/\rho^{n}_\edge$.
The factor $d$ should be suppressed to be consistent with what would be obtained by a finite volume discretization of this elliptic equation, if this latter was derived in the time semi-discrete setting: this fact is linked with the well-known non-consistency of the Rannacher-Turek or Crouzeix-Raviart discretization of the Darcy problem. 

\medskip
Then, equation \eqref{step31_3} is solved at the same time as the third equation of \eqref{proj-alg3} by a Newton's algorithm.
Once $p_{k+1}^{n+1}$ is known, the first relation of (\ref{proj-alg3}) gives the updated value of the velocity:
\begin{equation}
u_{k+1}^{n+1}=\tilde u^{n+1}-\dt\ (\matM_{\rho^{n}})^{-1}\,\matB^t\,(p_{k+1}^{n+1}-p^{n})
\label{step32_3}\end{equation}
As, to preserve the positivity of the density, we want to use in the mass balance the value of the density upwinded with respect to $u^{n+1}$, equations \eqref{step31_3} and \eqref{step32_3} are not decoupled, by contrast with what happens in usual projection methods.
They are thus solved in sequence, performing the upwinding with respect to $u_k^{n+1}$ in \eqref{step31_3} and then updating the velocity by \eqref{step32_3}, up to convergence.


\subsection{Spatial discretization of the correction step for $y$}\label{sec:y_correction}

To be consistent with the discretization of the first part of the gas mass balance, the correction step for $y$ is discretized by the finite volume method, and the resulting discrete problem reads: 
\begin{equation}
\begin{array}{l} \displaystyle
|K|\,\frac{\rho_{K}^{n+1}\,y_{K}^{n+1}-z_{K}^{n+1}}{\dt}
+\sum_{\edge=K|L}(G_{\edge,K}^{n+1})^+\,g(y_K^{n+1},y_L^{n+1})-(G_{\edge,K}^{n+1})^-\,g(y_L^{n+1},y_K^{n+1})
\hspace{10ex} \\[2ex] \displaystyle \hfill
+D\sum_{\edge=K|L}\frac{|\edge|}{d_\edge} (y_K^{n+1}-y_L^{n+1})=0
\end{array} 
\label{splitting_disc}
\end{equation}
In this relation, for all edge $\edge=K|L$, $d_\edge$ is the Euclidean distance between two points $x_K$ and $x_L$ of the adjacent meshes $K$ and $L$, supposed to be such that the segment $[x_K,\,x_L]$ is perpendicular to $K|L$.
These points may be defined as follows: if the control volume $K$ is a rectangle or a cuboid, $x_K$ is the barycenter of $K$; if the control volume $K$ is a simplex, $x_K$ is the circumcenter of the vertices of $K$.
Note that, in this latter case, the condition $x_K \in K$ implies some geometrical constraints for $K$.
Of course, in the cases where the diffusion coefficient $D=0$, these limitations are useless.

\medskip 
The quantities $(G_{\edge,K}^{n+1})^+$ and $(G_{\edge,K}^{n+1})^-$ are defined as $(G_{\edge,K}^{n+1})^+=\max(G_{\edge,K}^{n+1},0)$ and  $(G_{\edge,K}^{n+1})^-=-\min(G_{\edge,K}^{n+1},0)$ respectively, with $G_{\edge,K}^{n+1}$ given by:
\[
G_{\edge,K}^{n+1} = \rho_{\edge,{\rm up}}^{n+1} \ \int_\edge u_r^{n+1}  \cdot n_\edge
\]
where $\rho_{\edge,{\rm up}}^{n+1}$ stands for $\rho_{\edge,{\rm up}}^{n+1}=\rho_K^{n+1}$ if $u_\edge^{n+1}\cdot n_\edge \geq 0$ and $\rho_{\edge,{\rm up}}^{n+1}=\rho_L^{n+1}$ otherwise.
Note that this upwind choice with respect to $u^{n+1}$ has no theoretical justification: in fact, the developments of this paper hold with any discretization for this density, and we use here the same discretization as in the mass balance simply to make the informatic implementation easier.
The function $g(\cdot,\cdot)$ corresponds to an approximation of $\varphi(y)= \max [\,y\,(1-y),\ 0\,]$ by a monotone numerical flux function.
Let us recall the definition of this latter notion \cite{eym-00-fin}:
\begin{dfntn}[Monotone numerical flux function]\label{num_flux}
Let the function $g(\cdot,\cdot) \in C(\xR^2,\, \xR)$ satisfy the following assumptions:
\begin{enumerate}
\item $g(a_1,\,a_2)$ is non-decreasing with respect to $a_1$ and non-increasing with respect to $a_2$, for any real numbers $a_1$ and $a_2$,
\item $g(\cdot,\cdot)$ is Lipschitz continuous with respect to both variables over $\xR$,
\item $g(a_1,\,a_1)=\varphi(a_1)$, for any $a_1\in \xR$.
\end{enumerate} 
Then $g(\cdot,\cdot)$ is said to be a monotone numerical flux function for $\varphi(\cdot)$.
\end{dfntn}
Several choices are possible for the numerical flux function $g(\cdot,\cdot)$ and we refer to \cite{eym-00-fin} for some examples and references.
We adopt here the following simple flux-splitting formula:
\[
g(a_1,a_2)=g_1(a_1)+ g_2(a_2)
\]
where $g_1(a_1)=a_1$ if $a_1 \in [0,1]$ and $g_1(a_1)=0$ otherwise, and $g_2(a_2)=-(a_2)^2$ if $a_2 \in [0,1]$ and $g_2(a_2)=0$ otherwise.
Note that this choice does not exactly match the definition, as neither $g_1(\cdot)$ nor $g_2(\cdot)$ are continuous at $a_1=1$.
However, this is unimportant, as one can prove, even in this case, that the solution $y$ remains in the interval $(0,1]$, as stated in the following lemma which is a weaker version of the result proven in \cite[section 2]{gas-07-on}.

\begin{lmm}[Existence and uniqueness for a discrete solution]\label{exist_step_y}
Let us suppose that, $\forall K \in \mesh,\ \rho_K^{n+1} >0 $ and $z_K^{n+1}/\rho_K^{n+1} \in (0,1]$.
Then, there exists a unique solution to the considered discrete problem \eqref{splitting_disc}, and this solution verifies $y_K^{n+1} \in (0,1],\ \forall K\in \mesh$.
\end{lmm}


\subsection{Some properties of the scheme}

The following theorem gathers some properties of the scheme, which are essentially straightforward consequences of lemmas \ref{exist_step_v}, \ref{exist_step_p} and \ref{exist_step_y}.

\begin{thrm}[Properties of the scheme]
Let the density of the liquid phase be constant and the gas phase obeys the ideal gas law.
We suppose that the viscous term is dissipative (\ie\ $\forall v\in W_h,\ a_d(v,v)\geq 0$).
In addition, we assume that the initial density is positive and the initial gas mass fraction belongs to the interval $(0,1]$.
Then there exists a solution $(u^n)_{1\leq n \leq N}$, $(p^n)_{1\leq n \leq N}$, $(\rho^n)_{1\leq n \leq N}$, $(z^n)_{1\leq n \leq N}$ and $(y^n)_{1\leq n \leq N}$ to the scheme which enjoys the following properties, for all $n \leq N$:
\begin{itemize}
\item the unknowns lie in their physical range:
\[
\forall K \in \mesh, \qquad
\rho_K^n > 0, \qquad z_K^n >0, \qquad p_K^n >0, \qquad y_K^n \in (0,1]
\]
\item the total mass, the gas mass and, if $f_v=0$, the integral of the momentum are conserved:
\[
\begin{array}{l} \displaystyle
\sum_{K\in\mesh}|K|\ \rho^n_K = \sum_{K\in\mesh}|K|\ \rho^0_K
\\[3ex] \displaystyle
\sum_{K\in\mesh} |K|\ z^n_K = \sum_{K\in\mesh} |K|\ \rho^n_K y^n_K =\sum_{K\in\mesh} |K|\ \rho^0_K y^0_K
\\[3ex] \displaystyle
\sum_{\edge \in \edgesint} |D_\edge|\ \rho^{n-1}_\edge\ u_\edge^n= \sum_{\edge \in \edgesint} |D_\edge|\ \rho^{-1}_\edge\ u_\edge^0
\end{array}
\]
\end{itemize}
\end{thrm}

We now turn to another feature of the scheme, which, from numerical experiments, seems to be crucial for the robustness of the algorithm.
Let us suppose until the end of this section that the drift velocity $u_r$, the diffusive coefficient $D$ and the forcing term $f_v$ are set to zero.
In addition, we make abstraction of the boundary condition, \ie\ we momentarily reason as if the problem was posed in $\xR^n$.
Then the continuous problem enjoys the following property: if the initial velocity and the initial pressure are constant, let say $u=u_0$ and $p=p_0$ respectively, then they remain constant throughout the transient, while $\rho$ or $z$ are transported by this (constant) velocity; this solution corresponds to the transport of the contact discontinuity of the underlying hyperbolic system, the wave structure of which is quite similar to the Euler equations one \cite{gui-07-ada}.
The objective of the subsequent development is to prove that the numerical scheme considered in this paper presents the same behaviour: if, at the initial time, $u_K^0=u_0$ and $p_K^0=p_0$ for all $K\in\mesh$, then $p_K^{n+1}=p_0$ and $u_K^{n+1}=u_0$, for all $K\in\mesh$ and $n<N$.

\medskip
Let us assume that, at time $t=t_n$, the velocity $u^n$ and the pressure $p^n$ take the constant value $u_0$ and $p_0$ respectively.
We are now going to check that there exists a solution $u^{n+1}$, $p^{n+1}$, $z^{n+1}$ and $y^{n+1}$ to the scheme such that $u^{n+1}=u_0$ and $p^{n+1}=p_0$.
The discrete momentum balance equation reads, with a zero forcing term:
\[
\begin{array}{l} \displaystyle
\forall \edge \in \edgesint, \mbox{ for } 1 \leq i \leq d,
\\[2ex] \displaystyle \qquad
\frac{|D_\edge|} \dt\ (\rho_\edge^{n} \tilde u_{\edge,i}^{n+1}- \rho^{n-1}_\edge u_{\edge,i}^n)
+ \sum_{\stackrel{\scriptstyle \edged \in \edges(D_\edge),}{\scriptstyle \edged=D_\edge |  D_{\edge'}}}\frac 1 2\ \fluxd^{n}\ (\tilde u_{\edge,i}^{n+1} + \tilde u_{\edge',i}^{n+1})
-\int_{\Omega,h} p^{n}\ \nabla \cdot \varphi_{\edge}^{(i)} \,{\rm d}x
 = 0
\end{array}
\]
Replacing $u^{n}$ and $p^n$ by $u_0$ and $p_0$ respectively and taking $\tilde u^{n+1}_{\edge}=u_0$ for all $\edge\in\edgesint$, this system becomes:
\[
\forall \edge\in\edgesint, \ \edge=K|L,\qquad 
u_0 \ \left[\frac{|D_\edge|}{\dt} \, (\rho_\edge^{n}-\rho^{n-1}_\edge) 
+ \sum_{\scriptstyle \edged \in \edges(D_\edge)}\fluxd^{n}\right]=0
\]
which is verified thanks to the equivalence between mass balances over primal and dual meshes, as explained in section \ref{sec:momentum}. 
We now turn to the pressure correction step, which we recall:
\[
\left| \begin{array}{ll} \displaystyle
\frac{|D_\edge|}{\dt} \, \rho_\edge^{n}\,(u_\edge^{n+1} - \tilde u_\edge^{n+1}) 
+ |\edge| \, \left[ (p_K^{n+1}-p_K^{n})-(p_L^{n+1}-p_L^{n}) \right]\, n_{KL}= 0,
& \displaystyle
\forall \edge \in \edgesint,\ \edge=K|L
\\[2ex] \displaystyle
\frac{|K|}{\dt}\left[\varrho^{\,p,z}(p_K^{n+1},z_K^{n+1})-\rho^{n}_K\right]
\\[1ex] \displaystyle \hspace{5ex}
+\sum_{\edge \in \edges(K)} \left[(\vsp)^{n+1}
\, \varrho^{\,p,z}(p_K^{n+1},z_K^{n+1}) - (\vsm)^{n+1}\,\varrho^{\,p,z}(p_L^{n+1},z_L^{n+1})\right]=0,
\hspace{5ex} & \displaystyle
\forall K \in \mesh
\\[3ex] \displaystyle
\frac{|K|}{\dt}(z_K^{n+1}-z^{n}_K)+\sum_{\edge \in \edges(K)}\left[(\vsp)^{n+1}\, z^{n+1}_K - (\vsm)^{n+1}\,z^{n+1}_L\right]=0,
& \displaystyle
\forall K \in \mesh
\end{array}\right .
\]
Taking $u^{n+1}_{\edge}=u_0$ for all $\edge\in\edgesint$ and $p_K^{n+1}=p_0$ for all $K\in\mesh$, the left hand side of the first equation of this system vanishes.
Next, following \cite{gal-03-ont}, we remark that, at fixed pressure, the equation of state giving the density $\rho$ as a function of $z$ becomes an affine function:
\[
\rho=\varrho^{\,p,z}(p_0,z)=z\left(1-\frac{\rho_\ell\,a^2}{p_0}\right)+\rho_\ell
\]
Introducing this relation in the mass balance equation, we obtain:
\[\begin{array}{l} \displaystyle
\frac{|K|}{\dt}\left[\left(z_K^{n+1}-z_K^n\right)\,\left(1-\frac{\rho_\ell\,a^2}{p_0}\right)\right]
\\[3ex] \displaystyle \hspace{10ex}
+\sum_{\edge \in \edges(K)} (\vsp)^{n+1}\,\left[z_K^{n+1}\left(1-\frac{\rho_\ell\,a^2}{p_0}\right)+\rho_\ell\right]
-(\vsm)^{n+1}\, \left[z_L^{n+1}\left(1-\frac{\rho_\ell\,a^2}{p_0}\right)+\rho_\ell\right]=0
\end{array}
\]
which can be recast as:
\[
\left(1-\frac{\rho_\ell\,a^2}{p_0}\right)\,\left[\frac{|K|}{\dt}(z_K^{n+1}-z^{n}_K)+\sum_{\edge \in \edges(K)}\left[(\vsp)^{n+1}\, z^{n+1}_K - (\vsm)^{n+1}\,z^{n+1}_L\right]\right]
+ \rho_\ell \sum_{\edge \in \edges(K)} \vs^{n+1} =0
\]
which, as the last term vanishes for $u^{n+1}=u_0$, is exactly the same equation as the gas mass balance.
Thus, $u^{n+1}=u_0$, $p^{n+1}=p_0$, $z^{n+1}$ given by this latter equation and $y^{n+1}$ satisfying the correction step (which, for $u_r=0$ and $D=0$ becomes $\rho^{n+1} y^{n+1}=z^{n+1}$) is a solution to the scheme.
Consequently, provided that the solution is unique, the algorithm does preserve constant pressure and velocity through moving interfaces between phases, and transport this interface with this constant velocity.

\begin{rmrk}[More general boundary conditions]
The same property holds with a bounded computational domain when prescribing on the boundary either $u=u_0$ or a Neumann condition compatible with $u=u_0$ and $p=p_0$; this fact has been confirmed by numerical experiments, although we leave its proof beyond the scope of this presentation, to avoid the technicalities of the description of these latter discrete boundary conditions.
\end{rmrk}

\begin{rmrk}[On the choice of coupling the mass balance and the gas mass balance equations]
As in \cite{gas-07-on}, one may be tempted, specially for computing efficiency reasons, to use a fully fractional step algorithm, \ie\ to solve all the equations sequentially.
The central argument of the preceding development is that, with a fixed pressure, the quantity $\rho y$ is affine with respect to $\rho$, and this fact originates from the particular form of the equation of state.
Thus, for this argument to hold, it is mandatory for the density in the product $\rho y$ to be given by the equation of state $\rho=\varrho^{\,p,y} (y,p^\ast)$, where only the pressure may be taken at the previous time step or at the previous stage of the algorithm (indeed, when checking as below that the interface is transported with a fixed pressure, $p$ will be considered constant, in particular with respect to time).
Hence, the transport terms in the gas mass balance should read, in the time semi-discrete setting:
\[
\frac 1 \dt \ (\varrho^{\,p,y}(y^{n+1},p^n)\, y^{n+1} - \rho^n y^n) + \nabla \cdot \varrho^{\,p,y}(y^{n+1},p^n)\, y^{n+1} u^n + \dots =0
\]
But, in this case, the compatibility condition which yields a maximum principle for the advection operator, which here would read:
\[
\frac 1 \dt \ (\varrho^{\,p,y}(y^{n+1},p^n) - \rho^n) + \nabla \cdot \varrho^{\,p,y}(y^{n+1},p^n)\, u^n =0
\]
does not hold.
So it seems that an algorithm keeping $y$ within its physical bounds and transporting the interface at constant pressure and velocity necessarily couples the mixture and the gas mass balance.
\end{rmrk}


\section{The stability induced by the pressure forces work}\label{sec:p_work}

The aim of this section is to prove that the discretization at hand satisfies a stability bound which can be seen as the discrete analogue of the following equation :
\begin{equation}
-\int_\Omega p(x) \, \nabla \cdot u(x) \, {\rm d}x=\frac{d}{dt}\int_{\Omega} f(x) \dx,
\label{cont_pot}\end{equation}
where $f(x)$ stands for the volumetric free energy of the mixture.
The role played by this estimate in the theory which is developped here is twofold.
First, it provides an {\em a priori} bound for a class of discrete problems including the pressure correction step, which is the corner stone to prove the existence of a solution; this development is presented in appendix.
Second, it is crucial to derive stability results for the scheme.

\medskip
Throughout this section, we suppose that both the drift velocity and the gass fraction diffusion vanishes, so the overall and gas mass balance equations simply read:
\begin{eqnarray}
\frac{\partial\rho}{\partial t}+\dive(\rho u)=0 \label{ombal}
\\[1ex]
\frac{\partial z}{\partial t}+\dive(z u)=0 \label{gmbal}
\end{eqnarray}
Of course, stability results for the complete problem will {\em in fine} depend on the fact that the neglected terms in \eqref{gmbal} are dissipative with respect to the free energy; this point will be treated further.

\medskip
This section is organized as follows.
First, we prove this estimate in a general setting, \ie\ without specifying the equation of state for the fluid.
Then we explain as this theory applies to the case specifically adressed here, namely a constant density fluid and a gaseous phase obeying the ideal gas law.


\subsection{Abstract estimates}

The formal computation which allows to derive estimate \eqref{cont_pot} in the continuous setting is the following. 
The first assumption is that, through the (system of) equation(s) of state, the specific free energy can be expressed as a function of the mixture density and the gas partial density, which we write $f=f(\rho,z)$.
Then multiplying the mass balance equation by the derivative of $f$ with respect to $\rho$, the gas mass balance equation by the derivative of $f(\cdot,\cdot)$ with respect to $z$ and finally summing these relations, we obtain:
\[
\frac{\partial f}{\partial \rho}\ \left[\frac{\partial\rho}{\partial t}+\dive(\rho u)\right]
+ \frac{\partial f}{\partial z}\ \left[\frac{\partial z}{\partial t}+\dive(z u)\right]=0
\]
which yields:
\[
\frac{\partial}{\partial t} f(\rho(x,t),z(x,t))
+\frac{\partial f}{\partial \rho}\ \dive(\rho u) + \frac{\partial f}{\partial z}\ \dive(z u)=0
\]
Developping the divergence terms, we get:
\begin{equation}
\frac{\partial}{\partial t} f(\rho(x,t),z(x,t))
+ u \cdot \left[ \frac{\partial f}{\partial \rho}\ \grad \rho + \frac{\partial f}{\partial z}\ \grad z \right]
+ \dive u \left[\rho\ \frac{\partial f}{\partial \rho} + z\ \frac{\partial f}{\partial z} \right]=0
\label{preum}\end{equation}
The second term of this relation is equal to $u \cdot \grad f(\rho(x,t),z(x,t))$.
Adding and substracting $f\ \dive u$, we thus have:
\begin{equation}
\frac{\partial}{\partial t} f(\rho(x,t),z(x,t))+ \dive (f(\rho(x,t),z(x,t))\, u)
+ \dive u \left[\rho\ \frac{\partial f}{\partial \rho} + z\ \frac{\partial f}{\partial z} -f \right]=0
\label{deuz}\end{equation}
Since the integral of $\dive (f(\rho(x,t),z(x,t))\, u)$ over the computational domain vanishes thanks to the boundary conditions, this equation is the relation we are seeking, provided that the free energy is such that the following relation holds:
\[
\rho\ \frac{\partial f}{\partial \rho} + z\ \frac{\partial f}{\partial z} -f = p
\]

\medskip
We are going now to reproduce this computation at the discrete level.

%
%
\begin{thrm}[Stability due to the pressure work]\label{pot_el_pxz}
Let ${\cal C}$ be an open convex subset of $\xR^2$ and $f(\cdot,\cdot)$ be a convex continuously differentiable function from ${\cal C}$ to $\xR$.
We suppose that $(\rho_K)_{K\in\mesh}$, $(\rho_K^\ast)_{K\in\mesh}$, $(z_K)_{K\in\mesh}$ and $(z_K^\ast)_{K\in\mesh}$ are four families of real numbers such that, $\forall K \in \mesh$, $(\rho_K,z_K) \in {\cal C}$, $(\rho_K^\ast,z_K^\ast) \in {\cal C}$ and the following relations hold:
\begin{equation}
\left| \begin{array}{l} \displaystyle
\frac{|K|}{\dt}\,(\rho_K-\rho_K^\ast) +\sum_{\edge=K|L}\vs\,\rho_\edge=0
\\[2ex] \displaystyle
\frac{|K|}{\dt}\,(z_K-z_K^\ast)+\sum_{\edge=K|L} \vs \, z_\edge=0
\end{array}\right.
\label{balances}
\end{equation}
where $\rho_\edge$ and $z_\edge$ are given by $\rho_\edge=\rho_K$ and $z_\edge=z_K$ if $\vs\geq 0$, $\rho_\edge=\rho_L$ and $z_\edge=z_L$ otherwise.
Then the following estimate holds:
\[
\sum_{K\in\mesh} - p_K  \left[ \sum_{\edge=K|L} \vs \right] \geq
\sum_{K\in\mesh}  |K|\ \frac{f(\rho_K,z_K) - f(\rho_K^\ast,z_K^\ast)} \dt
\]
where the family of real numbers $(p_K)_{K\in\mesh}$ is given by:
\[
\forall K \in \mesh,\qquad 
p_K=\rho_K\ \frac{\partial f}{\partial \rho} (\rho_K,z_K)
+ z_K\ \frac{\partial f}{\partial z}(\rho_K,z_K) 
-f(\rho_K,z_K)
\]
\end{thrm}
\begin{proof}
Let us multiply the first relation of \eqref{balances} by the derivative with respect to $\rho$ of $f(\cdot,\cdot)$, the second one by the derivative with respect to $z$ of $f(\cdot,\cdot)$, both being evaluated at $(\rho_K,\, z_K)$, and sum:
\begin{equation}
\begin{array}{l} \displaystyle
\underbrace{
\frac{|K|}{\dt}\left[(\rho_K-\rho_K^\ast) \left(\frac{\partial f}{\partial\rho}\right)_{(\rho_K,\, z_K)}
                    +(z_K-z_K^\ast)\left(\frac{\partial f}{\partial z}\right)_{(\rho_K,\, z_K)}\right]
}_{T_{\partial/\partial t}}
\\ \displaystyle \hspace{20ex}
+\underbrace{
\left(\frac{\partial f}{\partial \rho}\right)_{(\rho_K,\, z_K)} \sum_{\edge=K|L} \vs \, \rho_\edge
+\left(\frac{\partial f}{\partial z}\right)_{(\rho_K,\, z_K)} \sum_{\edge=K|L} \vs \, z_\edge
}_{T_{div,K}}=0
\end{array}
\label{T1T2}\end{equation}
The second term of the previous relation, $T_{div,K}$, can be recast as:
\begin{equation}
\begin{array}{ll} \displaystyle
T_{div,K}
& \displaystyle
= \left(\frac{\partial f}{\partial \rho} \right)_{(\rho_K,\, z_K)}
\left[\sum_{\edge=K|L}\vs \, (\rho_\edge-\rho_K)+\rho_K\sum_{\edge=K|L} \vs \right]
\\[3ex] & \displaystyle \,
+ \left(\frac{\partial f}{\partial z}\right)_{(\rho_K,\, z_K)}
\left[\sum_{\edge=K|L}\vs \, (z_\edge-z_K)+z_K\sum_{\edge=K|L} \vs \right]
\end{array}
\label{deftdiv}\end{equation}
This relation is the discrete equivalent to equation \eqref{preum}: up to the multiplication by $1/|K|$, the first summations in the first term and the second term at the right hand side are the analogue of $u\cdot\grad \rho$ and $u\cdot\grad z$ respectively, while the second summations are the analogue of $\rho \dive u$ and $z\,\dive u$ respectively.
Adding and substracting $f(\rho_K,\, z_K)$, we obtain a discrete equivalent of relation \eqref{deuz}:
\[
\begin{array}{ll} \displaystyle
T_{div,K}
& \displaystyle
= \left(\frac{\partial f}{\partial \rho}\right)_{(\rho_K,\, z_K)} \sum_{\edge=K|L} \vs \, (\rho_\edge-\rho_K)
\\ & \displaystyle
+ \left(\frac{\partial f}{\partial z}\right)_{(\rho_K,\, z_K)} \sum_{\edge=K|L} \vs \, (z_\edge-z_K)
+ f(\rho_K,\, z_K) \sum_{\edge=K|L} \vs
\\ & \displaystyle
+ \left[ \rho_K\ \left(\frac{\partial f}{\partial \rho} \right)_{(\rho_K,\,z_K)}
   +z_K \left(\frac{\partial f}{\partial z} \right)_{(\rho_K,\, z_K)}
   -f(\rho_K,\, z_K) \right]
\sum_{\edge=K|L} \vs
\end{array}
\] 
In the last term, we recognize, as in the continuous setting, $p_K\sum_{\edge=K|L} \vs$.
The process will be completed if we put the first three terms of the right hand side in the divergence form. 
To this end, let us sum up the term $T_{div,K}$ over $K\in\mesh$ and reorder the summation:
\begin{equation}
\sum_{K\in\mesh}T_{div,K}=\sum_{K\in\mesh} p_K \left[ \sum_{\edge=K|L} \vs \right]
+\sum_{\edge\in\edgesint} T_{div,\edge}
\label{T3} \end{equation}
where, if $\edge=K|L$:
\[
\begin{array}{l} \displaystyle
T_{div,\edge}= \vs \left[ 
\left(\frac{\partial f}{\partial \rho}\right)_{(\rho_K,\, z_K)}(\rho_\edge-\rho_K)
+\left(\frac{\partial f}{\partial z}\right)_{(\rho_K,\, z_K)} (z_\edge-z_K)
+f(\rho_K,\, z_K)\right.
\\ \displaystyle \hspace{13ex}
\left.
-\left(\frac{\partial f}{\partial \rho}\right)_{(\rho_L,\, z_L)} (\rho_\edge-\rho_L)
-\left(\frac{\partial f}{\partial z}\right)_{(\rho_L,\, z_L)} (z_\edge-z_L)
-f(\rho_L,\, z_L) \right]
\end{array}
\]
In this relation, there are two possible choices for the orientation of $\edge$, \ie\ $K|L$ or $L|K$; we choose  this orientation in order to have $\vs \geq 0$.
The function $(\rho,\, z) \mapsto f(\rho,\, z)$ is by assumption continuously differentiable and convex on the convex set ${\cal C}$ containing both $(\rho_K,\, z_K)$ and $(\rho_L,\, z_L)$, so the technical lemma \ref{rhosigma_zsigma} hereafter applies and there exists $(\bar \rho_\edge,\, \bar z_\edge)$ in the segment $[(\rho_K,\, z_K),\ (\rho_L,\, z_L)]$ (itself included in ${\cal C}$) such that:
\begin{equation}
\left| \begin{array}{l}
\mbox{if }(\rho_K,\, z_K) \neq (\rho_L,\, z_L):
\\[2ex] \hspace{10ex}
\begin{array}{l} \displaystyle
 \left(\frac{\partial f}{\partial \rho}\right)_{(\rho_K,\, z_K)} (\bar\rho_\edge - \rho_K)
+\left(\frac{\partial f}{\partial z}   \right)_{(\rho_K,\, z_K)} (\bar z_\edge-z_K)
+ f(\rho_K,\, z_K)
\\[2ex] \displaystyle \hspace{15ex}
=\left(\frac{\partial f}{\partial \rho}\right)_{(\rho_L,\, z_L)} (\bar\rho_\edge-\rho_L)
+\left(\frac{\partial f}{\partial z}   \right)_{(\rho_L,\, z_L)} (\bar z_\edge-z_L)
+f(\rho_L,z_L)
\end{array}
\\[8ex]
\mbox{otherwise: }\quad
(\bar \rho_\edge,\, \bar z_\edge)=(\rho_K,\, z_K)=(\rho_L,\, z_L)
\end{array} \right.
\label{defrhob}\end{equation}
By definition, the choice $(\rho_\edge,z_\edge)=(\bar \rho_\edge,\bar z_\edge)$ is such that the term $T_{div,\edge}$ vanishes, which means that the first three terms at the right hand side of equation \eqref{deftdiv} are a conservative approximation of the quantity $\dive(f u)$ appearing in equation \eqref{deuz}, with the following expression for the flux:
\[
\begin{array}{l} \displaystyle
F_{\edge,K}= f_\edge \, \vs, \quad \mbox{with:}
\\[2ex] \displaystyle \hspace{15ex}
f_\edge= \left(\frac{\partial f}{\partial \rho} \right)_{(\rho_K,\, z_K)} (\bar\rho_\edge-\rho_K)
+\left(\frac{\partial f}{\partial z}\right)_{(\rho_K,\, z_K)} (\bar z_\edge-z_K)
+f(\rho_K,z_K)
\\[2ex] \displaystyle \hspace{18ex}
=\left(\frac{\partial f}{\partial \rho}\right)_{(\rho_L,\, z_L)} (\bar\rho_\edge-\rho_L)
+\left(\frac{\partial f}{\partial z}\right)_{(\rho_L,\, z_L)} (\bar z_\edge-z_L)
+f(\rho_L,z_L)
\end{array}
\]
Then the term $T_{div,\edge}$ can be rewritten as:
\[
T_{div,\edge}= 
\vs\ (\rho_\edge - \bar \rho_\edge)
\ \left[ \left(\frac{\partial f}{\partial \rho}\right)_{(\rho_K,\,z_K)} 
-\left(\frac{\partial f}{\partial \rho}\right)_{(\rho_L,\, z_L)} \right]
+\vs  \ (z_\edge - \bar z_\edge)
\ \left[ \left(\frac{\partial f}{\partial z}\right)_{(\rho_K,\, z_K)} 
-\left(\frac{\partial f}{\partial z}\right)_{(\rho_L,\, z_L)}\right]
\]
With the orientation taken for $\edge$, an upwind choice yields:
\[
T_{div,\edge}= 
\vs\ (\rho_K - \bar \rho_\edge)
\ \left[ \left(\frac{\partial f}{\partial \rho}\right)_{(\rho_K,\,z_K)} 
-\left(\frac{\partial f}{\partial \rho}\right)_{(\rho_L,\, z_L)} \right]
+\vs  \ (z_K - \bar z_\edge)
\ \left[ \left(\frac{\partial f}{\partial z}\right)_{(\rho_K,\, z_K)} 
-\left(\frac{\partial f}{\partial z}\right)_{(\rho_L,\, z_L)}\right]
\]
and, by the inequality of lemma \ref{rhosigma_zsigma} hereafter, $T_{div,\edge}$ can be seen to be non-negative.
Let us now turn to $T_{\partial/\partial t}$. 
As the function $(\rho,\, z) \mapsto f(\rho,\, z)$ is convex on the convex set ${\cal C}$ and both $(\rho_K,\, z_K)$ and $(\rho_K^\ast,\, z_K^\ast)$ belong to ${\cal C}$, we have:
\begin{equation}
T_{\partial/\partial t} \geq |K|\ \frac{f(\rho_K,\, z_K)- f(\rho_K^\ast, z_K^\ast)} \dt
\label{T1} \end{equation}
Then, summing for $K\in\mesh$ and using relations \eqref{T1T2}, \eqref{T3} and \eqref{T1} concludes the proof.
\end{proof}

%
%
In the course of the preceding proof, we used the following technical lemma.

\begin{lmm}\label{rhosigma_zsigma} 
Let ${\cal C}$ be an open convex subset of $\xR^2$, $f(\cdot,\cdot)$ be a convex continuously differentiable function from ${\cal C}$ to $\xR$ and $(\rho_1,\, z_1)$ and $(\rho_2,\, z_2)$ be two distinct elements of ${\cal C}$.
Then there exists $\zeta \in [0,1]$ such that $(\bar \rho, \, \bar z)= (1-\zeta)\ (\rho_1,\, z_1) + \zeta\ (\rho_2,\, z_2)$ satisfies the following relation:
\begin{equation}
\begin{array}{l}\displaystyle
f(\rho_1,z_1)
+\left(\frac{\partial f}{\partial \rho}\right)_{(\rho_1,\, z_1)} (\bar\rho-\rho_1)
+\left(\frac{\partial f}{\partial z}   \right)_{(\rho_1,\, z_1)} (\bar z-z_1)
=
\\[3ex] \displaystyle \hspace{20ex}
f(\rho_2,z_2)
+\left(\frac{\partial f}{\partial \rho}\right)_{(\rho_2,\, z_2)} (\bar\rho-\rho_2)
+\left(\frac{\partial f}{\partial z}   \right)_{(\rho_2,\, z_2)} (\bar z-z_2)
\end{array}
\label{defrhob_zb} \end{equation}
In addition, the following inequality holds:
\[
T=
(\rho_1-\bar \rho) \left[\left(\frac{\partial f}{\partial \rho}\right)_{(\rho_1,\, z_1)}
-\left(\frac{\partial f}{\partial \rho}\right)_{(\rho_2,\, z_2)} \right]
+
(z_1-\bar z) \left[\left(\frac{\partial f}{\partial z}\right)_{(\rho_1,\, z_1)}
-\left(\frac{\partial f}{\partial z}\right)_{(\rho_2,\, z_2)} \right] \geq 0
\]
\end{lmm}

\begin{proof}
Let us consider the function $g(\cdot)$ defined by:
\[
\zeta \mapsto f((1-\zeta)\ (\rho_1,\, z_1) + \zeta\ (\rho_2,\, z_2))
\]
By assumption, the function $g(\cdot)$ is defined over $[0,1]$, convex and continuously differentiable.
Moreover, it may be checked that equation \eqref{defrhob_zb} equivalently reads:
\[
g(0) + g'(0)\ \zeta = g(1) + g'(1)\ (\zeta-1)
\]
or, reordering terms:
\[
\left[g'(1)-g'(0)\right]\ \zeta = g(0)-(g(1)-g'(1))
\]
As $g(\cdot)$ is convex, if $g'(1)=g'(0)$, the function $g(\cdot)$ is affine and $g(0)-(g(1)-g'(1))$ vanishes, so the preceding relation is satisfied with any value of $\zeta$.
Otherwise, the preceding relation allows to compute $\zeta$ and, still by convexity of $g(\cdot)$, both $g'(1)-g'(0)$ and $g(0)-(g(1)-g'(1))$ is positive, and so is $\zeta$.
Still in this second case, this relation equivalently reads:
\[
\left[g'(1)-g'(0)\right]\ (\zeta-1)=g(0) + g'(0)-g(1)
\]
which, as $g(0) + g'(0)-g(1)$ is negative, shows that $\zeta \leq 1$.
Finally, the quantity $T$ simply reads $\zeta \left[ g'(1)-g'(0)\right]$, and is thus non-negative.
\end{proof}

\medskip
\begin{rmrk}[Discretization of the convective terms and conservation of the entropy]
From the above computation, it appears that the choice of $\bar \rho_\edge$ and $\bar z_\edge$ defined by equation \eqref{defrhob}, for the convective terms in the mass and the gas mass balance equations, is a convenient one to obtain an exact discrete counterpart of the continuous identity \eqref{cont_pot}, and thus, {\em in fine}, to build a scheme exactly conserving the entropy.
The upwind choice yields a dissipation, and nothing can be said for the centered one.
\end{rmrk}


\subsection{The case of a constant density liquid and an ideal gas}

Let us suppose that $\rho_\ell$ is constant and $\rho_g$ is linearly increasing with the pressure:
\[
\rho_g=\frac 1 {a^2}\ p
\]
where $a$ is a positive real number (from a physical point of view, it is the sound velocity in a pure gaseous isothermal flow).
For any positive $\rho$ and $z$ such that $z-\rho+\rho_\ell >0 $, the relation \eqref{eos} giving the mixture density as a function of the gas mass fraction and the phasic densities may be recast under the following form:
\begin{equation}
\frac 1 {a^2}\ p=\varrho^{\,\rho,z}_g (\rho,z)=\frac{z \ \rho_\ell}{z+\rho_\ell-\rho}
\label{rho_g}
\end{equation}
Let us define the volumetric free energy of the mixture by:
\begin{equation}
f(\rho,z)=a^2\ z\ \log(\varrho^{\,\rho,z}_g(\rho,z))
\label{f_rhoconst}\end{equation}
This function is continuously differentiable over the convex subset of $\xR^2$:
\begin{equation}
{\cal C}=\{(\rho,\, z) \in \xR^2 \mbox{ s.t. } \rho>0,\ z>0,\ z-\rho+\rho_\ell >0 \}
\label{def_C}\end{equation}
We are now going to show that it verifies the other two assumptions of theorem \ref{pot_el_pxz}, namely that $f(\cdot)$ is convex and satisfies the identity:
\[
T_p=\rho\ \frac{\partial f}{\partial \rho}+ z\ \frac{\partial f}{\partial z}-f =p
\]
This latter relation can be proven without referring to the specific form of $f(\cdot,\cdot)$, making use of the following property, which would be verified by the volumetric free energy function associated to any mixture composed of a constant density liquid phase and a barotropic gaseous phase:
\[
f(\rho,z)=z\ f_g(\varrho^{\,\rho,z}_g(\rho,z)) \qquad \mbox{with :} \quad f_g'(s)=\frac{\wp(s)}{s^2}
\]
where $\wp(\cdot)$ is the function giving the pressure as a function of the gas density (thus, in particular, $f'(\rho_g)=p/\rho_g^2$) and $f_g(\cdot)$ stands for the specific free energy of the gaseous phase.
Developping the derivatives and using the definition of $f_g(\cdot)$, we get:
\begin{equation}
T_{p}=
\rho\,z\, f_g'(\rho_g) \frac{\partial \varrho^{\,\rho,z}_g}{\partial \rho}
+ z^2 \,f_g'(\rho_g) \frac{\partial \varrho^{\,\rho,z}_g}{\partial z} 
+ z f_g(\rho_g) -z f_g(\rho_g)=
z \, \frac p {\rho_g^2} 
\left[ \rho\,\frac{\partial \varrho^{\,\rho,z}_g}{\partial \rho} + z \frac{\partial \varrho^{\,\rho,z}_g}{\partial z} \right]
\label{Tp}
\end{equation}
From the expression \eqref{rho_g}, we have:
\begin{equation}
\frac{\partial \varrho^{\,\rho,z}_g}{\partial \rho}=\frac{\rho_g^2}{\rho_\ell \ z}\qquad \mbox{ and } \quad
\frac{\partial \varrho^{\,\rho,z}_g}{\partial z}=\frac{\rho_g^2 (\rho_\ell-\rho)}{\rho_\ell \ z^2}
\label{d_zfg}
\end{equation}
Substituting in \eqref{Tp} leads to:
\[
T_p=\rho_g^2 f'_g(\rho_g)=p
\]
For proving the convexity of $f(\cdot,\cdot)$, we return to its explicit form:
\[
f(\rho,z)=a^2\ z\ \log \left(\frac{z \ \rho_\ell}{z+\rho_\ell-\rho}\right)
\] 
Differentiating twice this expression, we get:
\[
\frac{\partial^2 f}{\partial\rho^2}= a^2\,\frac{z}{(z+\rho_\ell+\rho)^2}, \qquad
\frac{\partial^2 f}{\partial z^2}  = a^2\,\frac{(\rho_\ell-\rho)^2}{z \, (z+\rho_\ell+\rho)^2},\qquad
\frac{\partial^2 f}{\partial\rho\partial z} = \frac{\partial^2 f}{\partial z \partial \rho}=
a^2\,\frac{\rho_\ell-\rho}{(z+\rho_\ell+\rho)^2}
\]
It is thus easy to check that the determinant of the Hessian matrix $A$ of $f(\cdot,\cdot)$ is zero while its trace is positive.
One eigenvalue of $A$ is thus zero and the second one is positive, and $f(\cdot,\cdot)$ is convex.

%
%

\section{Stability analysis}\label{sec:stab_z}

The aim of this section is to provide some results concerning the stability of (\ie\ the conservation of the entropy by) the scheme considered in this paper.
First (section \ref{sec:u_r=0}), in the case where both the drift velocity $u_r$ and the diffusion coefficient for the mass fraction of the dispersed phase $D$ vanish (\ie\ for the homogeneous model), we prove that the entropy (\ie\ the usual entropy associated to the homogeneous model) is conserved by the scheme, up to a step of renormalization of the pressure which is precisely stated.
Note that this step, which was implemented for monophasic flows in \cite{gal-07-an}, could be added in the present scheme; however, we have chosen not to consider it further than in this theoretical section, as, in practice, its beneficial effects were not clear.
Second (section \ref{sec:u_r=darcy}), we show that, as in the continuous case, if the drift velocity is proportional to the gradient of the pressure, the drift term is dissipative with respect to the same entropy; for this property to hold, a particular discretization of the drift term has to be implemented.

\medskip
In this section, we use the following discrete norm and semi-norm:
\begin{equation}
\begin{array}{ll} \D
\forall v \in W_h, \qquad
& \D
\normLdiscd{\rho}{v} = \sum_{\edge \in \edgesint} |D_\edge|\ \rho_\edge \, |v_\edge|^2
\\[2ex]
\forall q \in L_h, \qquad
& \D
\snormundiscd{\rho}{q}= \sum_{\edge \in \edgesint,\ \edge=K|L} \frac{1}{\rho_\edge}\ \frac{|\edge|^2}{|D_\edge|}\, (q_K-q_L)^2
\end{array}
\end{equation}
where $\rho=(\rho_\edge)_{\edge \in \edgesint}$ is a family of positive real numbers.
The function $\normLdiscd{\rho}{\cdot}$ defines a norm over $W_h$, and $\snormundisc{\rho}{\cdot}$ can be seen as a weighted version of the $H^1$ semi-norm classical in the finite volume context \cite{eym-00-fin}.
The following relation links this latter semi-norm to the problem at hand:
\begin{equation}
\forall q \in L_h, \qquad (\matB\,\matM_{\rho}^{-1}\, {\matB}^t\, q,q)=\snormundiscd{\rho}{q}
\label{lapl_p}
\end{equation}
where $\matB^t$, $\matB$ and $\matM_{\rho}$ are the discrete gradient, (opposite of the) divergence and mass matrix defined is section \ref{sec:p_correction}.
A proof of this equality can be found in \cite[section 3.4]{gal-07-an}.


\subsection{First case: $u_r=0$, $D=0$} \label{sec:u_r=0}

With a zero drift velocity and a zero diffusion coefficient, the numerical scheme at hand reads, in the time semi-discrete setting:
\begin{itemize}
\item[1 -] solve for $\tilde u^{n+1}$
\begin{equation}
\frac{\rho^n\ \tilde u^{n+1}-\rho^{n-1}\ u^n}{\dt}+\dive(\rho^n \ u^n \otimes\tilde u^{n+1})+\nabla p^n 
-\dive \tau(\tilde u^{n+1})= f_v^{n+1}
\label{qdm_0}
\end{equation}
\item[2 -] solve for $p^{n+1}$, $u^{n+1}$, $\rho^{n+1}$ and $z^{n+1}$
\begin{equation}
\left| \begin{array}{l} \displaystyle
\rho^n \ \frac{u^{n+1}-\tilde u^{n+1}}{\dt}+\nabla(p^{n+1}-p^n)=0
\\[2ex] \displaystyle
\frac{\varrho^{\,p,z}(p^{n+1},\ z^{n+1})-\rho^n}{\dt}+\dive(\varrho^{\,p,z}(p^{n+1},\ z^{n+1})\ u^{n+1})=0
\\[2ex] \displaystyle
\frac{z^{n+1}-\rho^n y^n}{\dt}+\dive(z^{n+1} \ u^n)=0
\\[2ex] \displaystyle
\rho^{n+1}=\varrho^{\,p,z}(p^{n+1},z^{n+1})
\end{array}\right.
\label{projection_0}\end{equation}
\item[3 -] solve for $y^{n+1}$
\begin{equation}
\rho^{n+1}y^{n+1}=z^{n+1}
\label{correction_y_0}\end{equation}
\end{itemize}

\smallskip
\begin{prpstn}[A partial stability result]
Let the density of the liquid phase be constant, the gas phase obeys the ideal gas law and $f(\rho,z)$ be the corresponding volumetric free energy of the mixture, defined by \eqref{f_rhoconst}.
We suppose that the viscous term is dissipative (\ie\ $\forall v\in W_h,\ a_d(v,v)\geq 0$).
In addition, we assume that the density $\rho^n$ is positive and the gas mass fraction $y^n$ belongs to the interval $(0,1]$.
Let $\tilde u^{n+1}$, $u^{n+1}$, $p^{n+1}$, $z^{n+1}$ and $\rho^{n+1}$ be a solution to equations \eqref{qdm_0}-\eqref{projection_0}, whith a zero forcing term.
Then the following bound holds:
\begin{equation}
\begin{array}{l} \displaystyle
\hspace{-2ex}\frac 1 2 \normLdiscd{\rho^{n}}{u^{n+1}} + \int_\Omega f(\rho^{n+1},z^{n+1}) \, {\rm d}x
+ \dt\, a_d(\tilde u^{n+1},\tilde u^{n+1})
+ \frac{\dt^2}{2} \snormundiscd{\rho^n}{p^{n+1}}
\hspace{15ex} \\ \displaystyle \hfill 
\leq \frac 1 2 \normLdiscd{\rho^{n-1}}{u^n} + \int_\Omega f(\rho^n,\rho^n y^n)\, {\rm d}x
+ \frac{\dt^2}{2} \snormundiscd{\rho^n}{p^n}
\end{array}
\label{stabres_part}\end{equation}
\label{part_stab}\end{prpstn}

\begin{proof}
Multiplying each equation of the first step of the scheme \eqref{qdm_0} by the corresponding unknown (\textit{i.e} the corresponding component of the velocity $\tilde u^{n+1}$ on the corresponding edge $\sigma$) and summing over the edges and the components yields, by virtue of theorem \ref{VF1_2}:
\begin{equation}
\frac{1}{2\, \dt} \normLdiscd{\rho^{n}}{\tilde u^{n+1}} - \frac{1}{2\, \dt} \normLdiscd{\rho^{n-1}}{u^n}
+ a_d(\tilde u^{n+1},\tilde u^{n+1}) - \int_{\Omega,h} p^{n} \nabla \cdot \tilde u^{n+1}\, {\rm d}x \leq 0
\label{stab1_case1}
\end{equation}
On the other hand, the first relation of system equation \eqref{projection_0} reads, in algebraic setting:
\[
\frac{1}{\dt} \matM_{\rho^{n}} \, (u^{n+1} - \tilde u^{n+1})+ \matB^t\, (p^{n+1} -p^{n})=0
\]
Reordering this relation and multiplying by $\matM_{ \rho^{n}}^{-1/2}$ (recall that $\matM_{\rho^{n}}$ is diagonal), we obtain:
\[
\frac{1}{\dt}\, \matM_{\rho^{n}}^{1/2} u^{n+1} 
+ \matM_{\rho^{n}}^{-1/2}\,\matB^t\, p^{n+1} =
\frac{1}{\dt}\, \matM_{\rho^{n}}^{1/2} \tilde u^{n+1} 
+ \matM_{\rho^{n}}^{-1/2}\,\matB^t\, p^{n}
\]
Squaring this relation gives:
\[
\begin{array}{l}\D 
\left(\frac{1}{\dt}\, \matM_{\rho^{n}}^{1/2} u^{n+1} 
+ \matM_{\rho^{n}}^{-1/2}\,\matB^t\, p^{n+1},
\ \frac{1}{\dt} \, \matM_{\rho^{n}}^{1/2} u^{n+1} 
+ \matM_{\rho^{n}}^{-1/2}\,\matB^t\, p^{n+1}\right) = 
\\ \D \hspace{20ex}
\left(\frac{1}{\dt}\, \matM_{\rho^{n}}^{1/2} \tilde u^{n+1} 
+ \matM_{\rho^{n}}^{-1/2}\,\matB^t\, p^{n},
\ \frac{1}{\dt}\, \matM_{\rho^{n}}^{1/2} \tilde u^{n+1} 
+ \matM_{\rho^{n}}^{-1/2}\,\matB^t\, p^{n}\right)
\end{array}
\]
which reads:
\[
\begin{array}{l} \D 
\frac{1}{\dt^{2}}\, \left( \matM_{\rho^{n}} u^{n+1}, \ u^{n+1}\right) 
+ \left(\matM_{\rho^{n}}^{-1}\,\matB^t\, p^{n+1}, \ \matB^t\, p^{n+1}\right)
+\frac{2}{\dt} \left(u^{n+1},\  \matB^t\, p^{n+1}\right)=
\\ \D \hspace{20ex}
\frac{1}{\dt^{2}}\, \left(\matM_{\rho^{n}} \tilde u^{n+1}, \ \tilde u^{n+1}\right) 
+\left(\matM_{\rho^{n}}^{-1}\,\matB^t\, p^{n}, \ \matB^t\, p^{n}\right)
+\frac{2}{\dt}\, \left(\tilde u^{n+1}, \ p^{n}\right)
\end{array}
\]
Multiplying by $\dt/2$, remarking that, $\forall v \in W_h,\ (\matM_{\rho^{n}} v, \ v)= \normLdiscd{\rho^{n}}{v}$ and that, thanks to relation \eqref{lapl_p}, $\forall q \in L_h,\ (\matM_{\rho^{n}}^{-1}\,\matB^t\, q, \ \matB^t\, q)=
 (\matB \, \matM_{\rho^{n}}^{-1}\,\matB^t\, q,q)=\snormundiscd{\rho^{n}}{q}$, we get: 
\begin{equation}
\begin{array}{l} \displaystyle
\frac{1}{2 \dt} \normLdiscd{\rho^{n}}{u^{n+1}}
+\frac \dt 2\, \snormundiscd{\rho^{n}}{p^{n+1}}
+(u^{n+1}, \matB^t\, p^{n+1})
\\ \hspace{20ex} \displaystyle
-\frac{1}{2 \dt} \normLdiscd{\rho^{n}}{\tilde u^{n+1}}
-\frac \dt2 \,\snormundiscd{\rho^{n}}{p^{n}}
- (\tilde u^{n+1}, \matB^t\, p^{n}) =0
\end{array}
\label{stab2_case1}\end{equation}
The quantity $-(\tilde u^{n+1}, \matB^t\, p^{n})$ is nothing more than the opposite of the term  $\D \int_{\Omega,h} p^{n} \nabla \cdot \tilde u^{n+1}\, {\rm d}x$ appearing in \eqref{stab1_case1}, so summing \eqref{stab1_case1} and \eqref{stab2_case1} makes these terms disappear, leading to:
\[
\begin{array}{l} \displaystyle
\frac{1}{2 \dt} \normLdiscd{\rho^{n}}{u^{n+1}}
- \frac{1}{2\, \dt} \normLdiscd{\rho^{n-1}}{u^n} 
+ a_d(\tilde u^{n+1},\tilde u^{n+1})
\\[2ex] \hspace{30ex} \displaystyle
+\frac \dt 2 \, \snormundiscd{\rho^{n}}{p^{n+1}}
-\frac \dt 2 \, \snormundiscd{\rho^{n}}{p^{n}}
+ (u^{n+1}, \matB^t\, p^{n+1})\leq 0 
\end{array}
\]
Finally, $(u^{n+1}, \matB^t\, p^{n+1})$ is precisely the pressure work which is likely to be bounded by the time derivative of the volumetric free energy of the mixture.
We know from theorem \ref{existence_pxz} that any solution to the system \eqref{projection_0} satisfies $\rho^{n+1}>0$, $z^{n+1}>0$ and $p^{n+1}>0$.
In view of the different forms of the equation of state gathered in \eqref{rhomixture}, this implies that this solution belongs to the convex set ${\cal C}$ defined by \eqref{def_C}, inside which the free energy is well defined, regular and convex.
Hence, with this solution, theorem \ref{pot_el_pxz} indeed applies and we get:
\[
\begin{array}{l} \displaystyle
\frac{1}{2 \dt} \normLdiscd{\rho^{n}}{u^{n+1}} 
+ a_d(\tilde u^{n+1},\tilde u^{n+1})
+\frac \dt 2\,\snormundiscd{\rho^n}{p^{n+1}}
+\frac{1}{\dt}  \int_\Omega f(\rho^{n+1},z^{n+1})\, {\rm d}x
\hspace{5ex} \\[2ex] \displaystyle \hfill
\leq \frac{1}{2\, \dt} \normLdiscd{\rho^{n-1}}{u^n} 
+\frac \dt 2\,\snormundiscd{\rho^n}{p^n} 
+\frac{1}{\dt} \int_\Omega f(\rho^n,\rho^n y^n)\, {\rm d}x 
\end{array}
\]
which concludes the proof.
\end{proof}


\begin{thrm}[Stability of the scheme, case $u_r={\cal D}=0$]
Let the density of the liquid phase be constant, the gas phase obeys the ideal gas law and $f(\rho,z)$ be the corresponding volumetric free energy of the mixture, defined by \eqref{f_rhoconst}.
We suppose that the viscous term is dissipative (\ie\ $\forall v\in W_h,\ a_d(v,v)\geq 0$).
In addition, we assume that the initial density is positive and the initial gas mass fraction belongs to the interval $(0,1]$.

\medskip
We now add to the scheme \eqref{qdm_0}-\eqref{correction_y_0} the following renormalization step of the pressure, to be performed at the very beginning of the time step, before the velocity prediction step:
\[
\mbox{Solve for } \tilde p^{n+1}: \qquad
-\nabla \cdot \left(\frac{1}{\rho ^n} \nabla \tilde p^{n+1}\right) = 
-\nabla \cdot \left(\frac{1}{\sqrt{\rho^n\,\rho^{n-1}}}\nabla p^n \right)
\]
or, in algebraic setting:
\[
\matB\,\matM^{-1}_{\rho^n}\, {\matB}^t \ \tilde p^{n+1} = 
\matB\,\matM^{-1}_{\sqrt{\rho^n\,\rho^{n-1}}}\, {\matB}^t \ p^n
\]
Accordingly, the pressure used in the velocity prediction step must be changed to $\tilde p^{n+1}$.

\medskip
Let $(\tilde u^n)_{0\leq n \leq N}$, $(u^n)_{0\leq n \leq N}$, $(p^n)_{0\leq n \leq N}$, $(z^n)_{0\leq n \leq N}$ and $(\rho^n)_{0\leq n \leq N}$ be the solution to this scheme, whith a zero forcing term.
Then the following entropy conservation result holds for $0 \leq n < N$:
\begin{equation}
\begin{array}{l} \displaystyle
\hspace{-2ex}\frac 1 2 \normLdiscd{\rho^{n}}{u^{n+1}} + \int_\Omega f(\rho^{n+1},z^{n+1}) \, {\rm d}x
+ \dt \sum_{k=1}^{n+1} a_d(\tilde u^k,\tilde u^k)
 + \frac{\dt^2}{2} \snormundiscd{\rho^{n}}{p^{n+1}}
\hspace{25ex} \\ \displaystyle \hfill
\leq \frac 1 2 \normLdiscd{\rho^0}{u^0} + \int_\Omega  z^{0} \, f_g(\rho^{0},z^{0})\, {\rm d}x
+ \frac{\dt^2}{2} \snormundiscd{\rho^0}{p^0}
\end{array}
\label{stabres_case1}
\end{equation}
\end{thrm}

\begin{proof}
By the same proof as for the scheme without the pressure renormalization step, we get:
\[
\begin{array}{l} \displaystyle
\frac{1}{2 \dt} \normLdiscd{\rho^{n}}{u^{n+1}} 
+ a_d(\tilde u^{n+1},\tilde u^{n+1})
+\frac \dt 2\,\snormundiscd{\rho^n}{p^{n+1}}
+\frac{1}{\dt}  \int_\Omega f(\rho^{n+1},z^{n+1})\, {\rm d}x
\hspace{15ex} \\[2ex] \displaystyle \hfill
\leq \frac{1}{2\, \dt} \normLdiscd{\rho^{n-1}}{u^n} 
+\frac \dt 2\,\snormundiscd{\rho^n}{\tilde p^{n+1}} 
+\frac{1}{\dt} \int_\Omega f(\rho^n,\rho^n y^n)\, {\rm d}x 
\end{array}
\]
and the conclusion follows by summing over the time steps, remarking that $z^{n+1}=\rho^{n+1}y^{n+1}$ and, thanks to the renormalization step (see \cite{gal-07-an} for a detailed computation):
\[
\snormundiscd{\rho^n}{\tilde p^{n+1}} \leq \snormundiscd{\rho^{n-1}}{p^n}
\]
\end{proof}

Note that a similar pressure renormalization step has already been introduced for variable density incompressible flows \cite{gue-00-proj}.


\subsection{Dissipativity of the drift term} \label{sec:u_r=darcy}

We address in this section the case where the drift velocity is given by the Darcy-like closure relation \eqref{u_r}:
\[
u_r =  \frac 1 \lambda\ (1-\alpha_g)\,\alpha_g\,\frac{\varrho_g(p)-\rho_\ell}{\rho}\ \grad p
\]
In this relation, $\lambda$ is a positive phenomenological coefficient and $\alpha_g$ is the void fraction, which can be expressed as a function of the unknowns used in the scheme as $\alpha_g=z/\varrho_g(p)$.
We recall the spatial discretization of the drift term in the correction step for the gas mass fraction $y$, namely $\dive(\rho\ y\,(1-y)\ u_r)$, given in section \ref{sec:y_correction}:
\[
\sum_{\edge=K|L}G_{\edge,K}^+\,g(y_K,y_L)-G_{\edge,K}^-\,g(y_L,y_K)
\]
where the function $g(\cdot,\cdot)$ corresponds to an approximation of $\varphi(y)= \max [\, y\,(1-y),\ 0\,]$ by a monotone numerical flux function, $G_{\edge,K}^+=\max(G_{\edge,K},0)$, $G_{\edge,K}^-=-\min(G_{\edge,K},0)$ and $G_{\edge,K}$ is an approximation for the flux of $\rho \ u_r$ through the edge $\edge=K|L$.
With the closure relation \eqref{u_r} for $u_r$, a natural discretization for this quantity reads:
\begin{equation}
G_{\edge,K} = |\edge|\ \rho_{\edge,{\rm up}}
\left[\frac{\alpha_g\,(1-\alpha_g)}{\lambda} \, (\rho_g-\rho_\ell)\right]_\edge 
(p_K-p_L)
\label{G_grad_p}\end{equation}
where $\rho_{\edge,{\rm up}}$ is a density on $\edge$, for which, for pratical implementation reasons, we choose an upwind approximation with respect to the mean velocity $u$.
The goal of this section is to show that it is possible to approximate:
\[
\left[\frac{\alpha_g\,(1-\alpha_g)}{\lambda} \, (\rho_g-\rho_\ell)\right]_\edge 
\]
in such a way that this drift term is dissipative with respect to the entropy of the system.

\medskip
We begin this section by stating a consequence of the equation of state for the mixture which is central to the present development.

\begin{lmm}
Let the density of the liquid phase be constant, the gas phase obeys the ideal gas law, $f(\rho,z)$ be the corresponding volumetric free energy of the mixture, defined by \eqref{f_rhoconst}, and $h(\rho,z)$ be the partial derivative of $f(\cdot,\cdot)$ with respect to the second variable $z$.
Then the following results hold:
\begin{enumerate}
\item $h(\cdot,\cdot)$ only depends on the pressure, \ie\ there exists a function $h_p(\cdot)$ such that, for $\rho$ and $z$ in the convex set ${\cal C}$ defined by \eqref{def_C},  $h(\rho,z)=h_p(\wp(\rho,z))$, where $\wp(\cdot,\cdot)$ is the function giving the pressure as a function of $\rho$ and $z$:
\[
p=\wp(\rho,z)=a^2\, \varrho^{\,\rho,z}_g (\rho,z)=a^2\, \frac{z \ \rho_\ell}{z+\rho_\ell-\rho}
\]
\item the derivative of $h_p(\cdot)$ is given by:
\[
h'_p(p)=\frac{\rho_\ell-\rho_g(p)}{\rho_\ell\,\rho_g(p)}
\]
\item for any positive real numbers $p_1$ and $p_2$ such that $p_1 < p_2$, there exists $p_{1,2} \in [p_1,p_2]$ such that:
\[
h'_p(p_{1,2})\ \frac{h_p(p_1)-h_p(p_2)}{p_1-p_2}\geq 0
\]
\end{enumerate}
\label{choice_edge}\end{lmm}

\begin{proof}
As $(\rho,z)\in {\cal C}$, the pressure or, equivalently, the gas density $\rho_g$ can be expressed as a function of $(\rho,z)$ by $\rho_g=\varrho^{\,\rho,z}_g(\rho,z)$.
By the definition of $f(\cdot,\cdot)$, we thus have:
\[
h(\rho,z)=\frac{\partial f}{\partial z}=
f_g(\rho_g)+z\ \frac{\partial f_g}{\partial z}=
a^2\,log\left(\frac{p}{a^2}\right)+ z\frac{\partial f_g}{\partial \rho_g} \frac{\partial \varrho^{\,\rho,z}_g}{\partial z}
 \]
Then using the expression \eqref{d_zfg} of the derivative of $\varrho^{\,\rho,z}_g(\cdot,\cdot)$ with respect to the second variable, we get:
\[ 
h(\rho,z)=a^2\,log \left(\frac{p}{a^2}\right)+ p\,\frac{\rho_\ell -\rho}{\rho_\ell\, z}
 \]
Using the fact that $\rho=(1-\alpha_g)\rho_\ell+\alpha_g\,\rho_g$ and thus $\displaystyle \rho_\ell -\rho=\alpha_g\,(\rho_\ell -\rho_g)=\frac z {\rho_g}\ (\rho_\ell -\rho_g)$, we have:
\[ 
 h(\rho,z)=a^2\,log\left(\frac{p}{a^2}\right)+ p\,\frac{\rho_\ell -\rho_g}{\rho_\ell\, \rho_g}
 \] 
By definition of $\rho_g$, \ie\ $\rho_g=p/a^2$, we thus get: 
\[ 
 h(\rho,z)=a^2 \left[ log\left(\frac{p}{a^2}\right)+ \frac{\rho_\ell -p/a^2}{\rho_\ell}\right]=h_p(p)
 \]
Taking the derivative of this relation yields the desired expression for $h'_p(\cdot)$ and, as $h_p(\cdot)$ is continuously differentiable in $[p_1,p_2]$, the existence of $p_{1,2}$ follows by Lagrange's theorem.
\end{proof}

\medskip
We are now in position to state and prove the following stability result.
\begin{prpstn}
Let the density of the liquid phase be constant, the gas phase obeys the ideal gas law and $f(\rho,z)$ be the corresponding volumetric free energy of the mixture, defined by \eqref{f_rhoconst}.
Let ${\cal C}$ be the convex set defined by \eqref{def_C} and $(\rho_K)_{K\in\mesh}$, $(y_K)_{K\in\mesh}$ and $(z_K)_{K\in\mesh}$ be such that, $\forall K \in \mesh$, $(\rho_K, \rho_K y_K) \in {\cal C}$, $(\rho_K, z_K) \in {\cal C}$, and the following relation is satisfied:
\begin{equation}
\frac{|K|}{\dt}\, (\rho_K \,y_K-z_K)
+\sum_{\edge=K|L} G_{\edge,K}^+\,g(y_K,y_L)-G_{\edge,K}^-\,g(y_L,y_K)=0
\label{balance_Y_gradp}\end{equation}
where $g(\cdot,\cdot)$ corresponds to an approximation of $\varphi(y)= \max[\, y\,(1-y),\ 0\,]$ by a monotone numerical flux function,
$G_{\edge,K}^+=\max(G_{\edge,K},0)$, $G_{\edge,K}^-=-\min(G_{\edge,K},0)$ and $G_{\edge,K}$ is given by the relation \eqref{G_grad_p}.
Then, if $g(y_K,y_L)\geq 0$ for all $\edge\in\edgesint$, $\edge=K|L$, there exists a discretization for the term:
\[
\left[\frac{\alpha_g\,(1-\alpha_g)}{\lambda} \, (\rho_g-\rho_\ell)\right]_\edge 
\]
in \eqref{G_grad_p} such that the following stability estimate holds:
\[
\frac{1}{\dt}\sum_{K\in\mesh}  |K|\ \left[f(\rho_K,\rho_K\,y_K) - f(\rho_K,z_K) \right]\leq 0
\]
which means that the drift term is dissipative with respect to the entropy of the system.
\label{potel_case2}\end{prpstn}

\begin{proof}
We multiply equation \eqref{balance_Y_gradp} by the partial derivative of $f(\cdot,\cdot)$ with respect to the second variable, taken at the point $(\rho_K,\rho_K \,y_K)$, and sum up over the control volumes of the mesh:
\[
\sum_{K\in\mesh} h(\rho_K, \rho_K \, y_K) \left[\frac{|K|}{\dt}\,(\rho_K \, y_K- z_K)
+\sum_{\edge=K|L}G_{\edge,K}^+\,g(y_K,y_L)-G_{\edge,K}^-\,g(y_L,y_K)\right]=T_1+T_2=0
\]
where $T_1$ and $T_2$ reads:
\[
\begin{array}{l}\displaystyle
T_1 = \sum_{K\in\mesh}\frac{|K|}{\dt}\ h(\rho_K,\rho_K\,y_K)\, \left[\rho_K \, y_K- z_K \right]
\\ \displaystyle
T_2 = \sum_{K\in\mesh}h(\rho_K,\rho_K\,y_K)
\left[\sum_{\edge=K|L}G_{\edge,K}^+\,g(y_K,y_L)-G_{\edge,K}^-\,g(y_L,y_K)\right]
\end{array}
\]
As the fonction $f(\cdot,\cdot)$ is convex, we have:
\begin{equation}
T_1\geq \frac{1}{\dt}\sum_{K\in\mesh}  |K|\ \left[ f(\rho_K,\rho_K\,y_K) - f(\rho_K,z_K)\right]
\label{t1}
\end{equation}
Let us turn to $T_2$. 
Reordering the sum, we get:
\[
T_2= \sum_{\sigma\in\edgesint}
|\edge|\ \rho_{\edge,{\rm up}} 
\, g_{\rm up}(y_K,y_L, u_r)
\left[\frac{\alpha_g\,(1-\alpha_g)}{\lambda} \, (\rho_g-\rho_\ell)\right]_\edge (p_K-p_L)
\left[h(\rho_K,\rho_K\,y_K)-h(\rho_L,\rho_L\,y_L)\right]
\]
where $g_{\rm up}(y_K,y_L, u_r)=g(y_K,y_L)$ if $u_r \geq 0$ and $g_{\rm up}(y_K,y_L, u_r)=g(y_L,y_K)$ if $u_r \leq 0$; in any case, we have, by assumption, $g_{\rm up}(y_K,y_L, u_r) \geq 0$.
We now choose, for the approximation of the quantity defined on $\edge$ in the preceding relation, an expression of the form:
\[
\left[\frac{\alpha_g\,(1-\alpha_g)}{\lambda} \, (\rho_g-\rho_\ell)\right]_\edge \equiv
\frac{(\alpha_g)_\edge\,\left[1-(\alpha_g)_\edge\right]}{\lambda} \, (\varrho_g(p_\edge)-\rho_\ell)
\]
where $(\alpha_g)_\edge$ stands for an approximation of the void fraction on $\edge$ which only needs here to be supposed non-negative.
Applying lemma \ref{choice_edge}, $T_2$ reads: 
\[
T_2= \sum_{\sigma\in\edgesint}
|\edge|\ \rho_{\edge,{\rm up}} 
\ g_{\rm up}(y_K,y_L, u_r)
\ \frac{(\alpha_g)_\edge\,(1-\alpha_g)_\edge}{\lambda} \ \rho_\ell\ \varrho_g(p_\edge)\ h'_p(p_\edge)\ (p_K-p_L) \ [h_p(p_K)-h_p(p_L)]
\]
If $p_K=p_L$, the term associated to $K|L$ in this sum vanishes.
Otherwise, from the third assertion of lemma \ref{choice_edge}, there exists $p_\edge \in [\min(p_K,p_L),\ \max(p_K,p_L)]$ such that the product $h'_p(p_\edge)\ (p_K-p_L) \ [h_p(p_K)-h_p(p_L)]$ is positive.
Since we choose $(\alpha_g)_\edge$ such that $(\alpha_g)_\edge \geq 0$, all the other quantities are positive, and this concludes the proof.
\end{proof}

The following proposition extends the stability result of the preceding section to the case $u_r\neq 0$.

\begin{prpstn}[Stability of the scheme, case $u_r\neq 0$]
Let the density of the liquid phase be constant, the gas phase obeys the ideal gas law and $f(\rho,z)$ be the corresponding volumetric free energy of the mixture, defined by \eqref{f_rhoconst}.
We suppose that the viscous term is dissipative (\ie\ $\forall v\in W_h,\ a_d(v,v)\geq 0$).
In addition, we assume that the density $\rho^n$ is positive and the gas mass fraction $y^n$ belongs to the interval $(0,1]$.
Let $\tilde u^{n+1}$, $u^{n+1}$, $p^{n+1}$, $z^{n+1}$, $\rho^{n+1}$ and $y^{n+1}$ be a solution to the equations of one time step of the scheme, whith a zero forcing term.
We suppose that the drift velocity is given by a the Darcy-like relation \eqref{u_r} and that the discretization of the correction step for the gas mass fraction $y^{n+1}$ is such that the stability result of proposition \ref{potel_case2} applies.
Then the following bound holds:
\[
\begin{array}{l} \displaystyle
\hspace{-2ex}\frac 1 2 \normLdiscd{\rho^{n}}{u^{n+1}} + \int_\Omega f(\rho^{n+1},\rho^{n+1} y^{n+1}) \, {\rm d}x
+ \dt \, a_d(\tilde u^{n+1},\tilde u^{n+1})
+ \frac{\dt^2}{2} \snormundiscd{\rho^n}{p^{n+1}}
\hspace{15ex} \\ \displaystyle \hfill 
\leq \frac 1 2 \normLdiscd{\rho^{n-1}}{u^n} + \int_\Omega f(\rho^n,\rho^n y^n)\, {\rm d}x
+ \frac{\dt^2}{2} \snormundiscd{\rho^n}{p^n}
\end{array}
\]
\label{stab_case2}\end{prpstn}

\begin{proof}
Proposition \ref{potel_case2} yields:
\[
\frac{1}{\dt}\sum_{K\in\mesh}  |K|\ \left[f(\rho_K^{n+1},\rho_K^{n+1}\,y_K^{n+1}) - f(\rho_K^n,z_K^n) \right]\leq 0
\]
The conclusion thus follows by summing this relation with the estimate of proposition \ref{part_stab}.
\end{proof}

Finally, note that, as in the preceding section, this partial stability result yields the same entropy decrease estimate for the whole scheme as in the preceding section if a renormalization step for the pressure is added to the scheme.

\begin{rmrk}[On the choice of the monotone numerical flux function]
As stated in section \ref{sec:algo_z}, we have adopted for the numerical tests presented hereafter the following flux-splitting formula:
\[
g(a_1,a_2)=g_1(a_1)+ g_2(a_2)
\]
where $g_1(a_1)=a_1$ if $a_1 \in [0,1]$ and zero otherwise and and $g_2(a_2)=-(a_2)^2$ if $a_2 \in [0,1]$ and zero otherwise.
This numerical monotone flux does not satisfy the hypothesis of proposition \ref{stab_case2}, as it is not always non-negative.
However, several other choices are possible for the numerical flux function $g(\cdot,\cdot)$ (\eg \cite{eym-00-fin}), and some of them solve this problem. 
Thanks to the fact that $\varphi(s)=s\,(1-s)$ is positive $\forall s \in[0,1]$, it is the case, for example, for the flux obtained with a one-dimensional Godunov scheme for each interface:
\[
g(a_1,a_2)=\left| \begin{array}{ll}
\D
\max\{\varphi(s), \,a_2\leq s \leq a_1\} & \mbox{if }a_2\leq a_1 \\
\D
\min\{\varphi(s), \,a_1\leq s \leq a_2\} & \mbox{if }a_1\leq a_2
                  \end{array}\right.
\]
\end{rmrk}


\section{Numerical results}\label{sec:numerical_z}

This section is devoted to numerical tests of the proposed scheme.
We first adress a problem built in such a way that it admits an analytical solution, to assess the convergence properties of the scheme.
Then several additional tests are performed, to check the stability of the algorithm and the quality of the results.


\subsection{Assessing the convergence against an analytic solution}

We address here a problem built by the so-called technique of manufactured solutions: the computational domain and the solution are chosen {\em a priori} and the initial conditions, the boundary conditions and the forcing terms are adjusted consequently. 
Let thus the computational domain be $\Omega=(0,1)\times (- 1/2,1/2)$, and the density and the momentum take the following expressions:
\[
\rho=1+\frac{1}{4}\,\sin(\pi t)\,\left[\cos(\pi x_1)-\sin(\pi x_2)\right]
\hspace{15ex}
\rho\,u= -\frac{1}{4} \cos(\pi t)\left[\begin{array}{l} \sin(\pi x_1) \\ \cos(\pi x_2) \end{array}\right]
\]
The pressure and the partial gas density are linked to the density by the equation of state \eqref{eos_z}, where the liquid density $\rho_{\ell}$ is set at $\rho_{\ell}= 5$ and the quantity $a^2$ in the equation of state of the gas \eqref{gas_eos} is given by $a^2=1$ (so $\rho_g=p$).
We choose the following expression for the unknowns $y$ and $z$:
\[
y=\frac{2.5 -0.5\ \rho}{4.5 \ \rho}
\hspace{15ex}
z=\rho\,y=\frac{2.5 -0.5\ \rho}{4.5}
\]
The relative velocity is constant and given by $u_r=(0,1)^t$ and the diffusive coefficient $D$ is set to $D=0.1$.
The analytical expression for the pressure is obtained from the equation of state (\ie\ relation \eqref{rho_g}).
These functions satisfy the mass balance equation; for the gas mass fraction and momentum balance, we add the corresponding right-hand side.
In this latter equation, we suppose that the divergence of the stress tensor is given by:
\[
\nabla \cdot \tau(u)= \mu\, \Delta u + \frac \mu 3\, \nabla\, \nabla \cdot u,
\qquad \mu=10^{-2}
\]
and we use for the viscous term the corresponding form for the bilinear form $a_d(\cdot,\cdot)$ (see section \ref{sec:momentum}).

\begin{figure}[htb]
\begin{center} \scalebox{0.4}{\includegraphics*{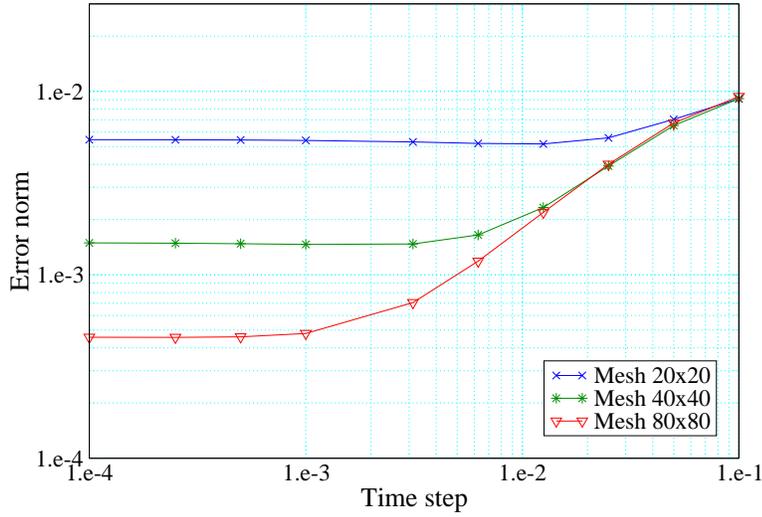}} \end{center}
\caption{Error for the velocity at $t=0.5$, as a function of the time step ($L^2$ norm).
\label{err_v}}
\end{figure}

\begin{figure}[htb]
\begin{center} \scalebox{0.4}{\includegraphics*{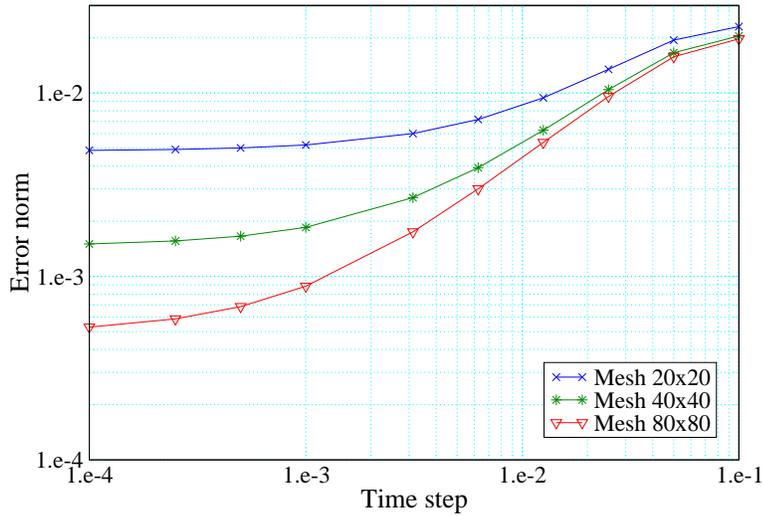}} \end{center}
\caption{Error for the pressure at $t=0.5$, as a function of the time step (discrete $L^2$ norm).
\label{err_p}}
\end{figure}

\begin{figure}[htb]
\begin{center} \scalebox{0.4}{\includegraphics*{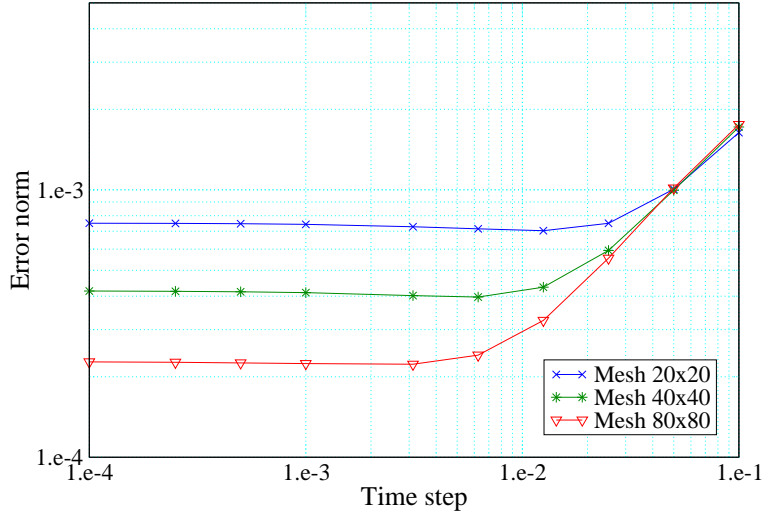}} \end{center}
\caption{Error for the gas mass fraction at $t=0.5$, as a function of the time step (discrete $L^2$ norm).
\label{err_y}}
\end{figure}

\medskip
Errors for the velocity, pressure and gas mass fraction obtained at $t=0.5$, as a function of the time step and for various meshings, are drawn on figure \ref{err_v}, figure \ref{err_p} and figure \ref{err_y}, respectively.
These errors are evaluated in the $L^2$ norm for the velocity and in the discrete $L^2$ norms for the pressure and the gas mass fraction.
Computations are made with $20 \times 20$, $40 \times 40$ and $80 \times 80$ uniform meshes (so with square cells and the Rannacher-Turek element).
For large time steps, these curves show a decrease corresponding to approximately a first order convergence in time, until a plateau is reached, due to the fact that errors are bounded by below by the residual spatial discretization error.
The value of the errors on this plateau then show a spatial convergence order close to one, which is consistent with the choice of an upwind discretization for the advection terms in the mass and gas mass fraction balance equations.

\subsection{Two-dimensional sloshing in cavity}

Two layers of non-miscible fluids (air and water) are superimposed with the lighter one on top of the heavier one.
The gravity (with $g=9.81 \, m.s^{-2}$) is acting in the vertical downward direction.
The length of the rectangular cavity is $L=1\,m$, the height of each layer is respectively $h_{\ell}=1\,m $ and  $h_g=1.25\,m $, so the total height of the box is $2.25\,m$. 
The water and air densities are respectively $\rho_{\ell}=1000 \, kg.m^{-3}$ and $\rho_g = p/a^2$ where $a^2$ is such that $\rho_g = 1.2 \, kg.m^{-3}$ at $p=10^5 \, Pa$. 
The diffusion coefficient $D$ and the drift velocity are set to zero.
A perfect slip condition is imposed on the whole boundary.
At initial time, both fluids are at rest, then the cavity is submitted to an horizontal acceleration given by $a_0 = 0.1\, m.s^{-2}$.

\medskip
In the case where both fluids are supposed incompressible and the convection and diffusion terms may be neglected, an analytical solution for the flow in a rectangular cavity is provided in \cite{cha-04-mod}.
In particular, the shape of the interface is given by the following relation:
\[
\xi=\frac{a_0}{g}\,\left[x-\frac{L}{2} + \sum_{n\geq 0}\frac{4}{L\, k^2_{2n+1}}\,\cos(\omega_{2n+1}\,t)\,\cos(k_{2n+1}\,t)\right]
\]
where the wave number $k_n$ is defined by:
\[
k_n=\frac{2\,\pi\,n}{L}
\]
and $\omega_n$ is given by:
\[
\omega_n^2=\frac{g\,k_n\,(\rho_{\ell}-\rho_g)}{\rho_g\,\coth(k_n\,h_g)+\rho_{\ell}\,\coth(k_n\,h_{\ell})}
\]
In practice, to compute this analytical solution, we perform the summation up to $n=200$.

\medskip
As, to remain in the domain of validity of the solution, the amplitude of the fluid oscillations must be very small, a very fine mesh is necessary near the free surface, to capture its motion.
The mesh is thus made of about $41\,000$ rectangular cells (with the Rannacher-Turek element) and, in the vertical direction, the space step is adapted in such a way that it is smaller near the interface between the two phases and equal to $\delta x_2= 0.0005\, m$, and increases when moving away the free surface, up to $\delta x_2 = 0.05\, m$ at the top and bottom sections.
In the horizontal direction, the mesh is uniform with step size $\delta x_1=1/70\, m$.
Calculations with different viscosities have been performed, these latter being supposed to vary with the mixture density: $\mu= \rho/100$, $\mu= \rho/1000$, $\mu= \rho/10000$.

\begin{figure}[htb]
\begin{center} \scalebox{0.6}{\includegraphics*{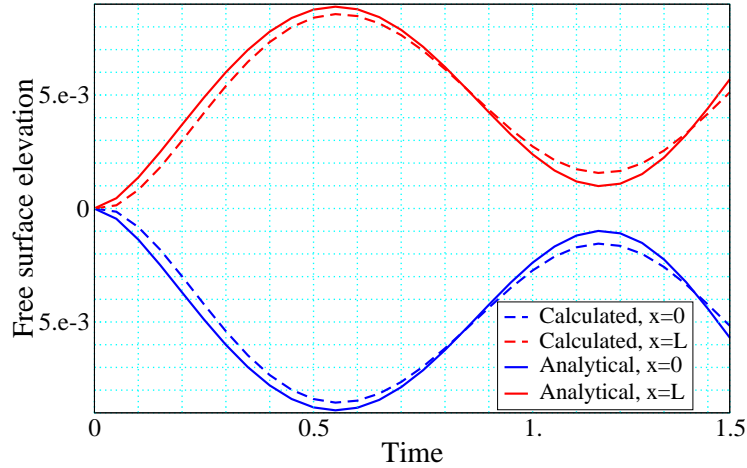}} \end{center}
\caption{Sloshing in cavity: analytical solution and numerical solution with $\mu= \rho/100$.
\label{ball_num_rho100}}
\end{figure}

\begin{figure}[htb]
\begin{center} \scalebox{0.6}{\includegraphics*{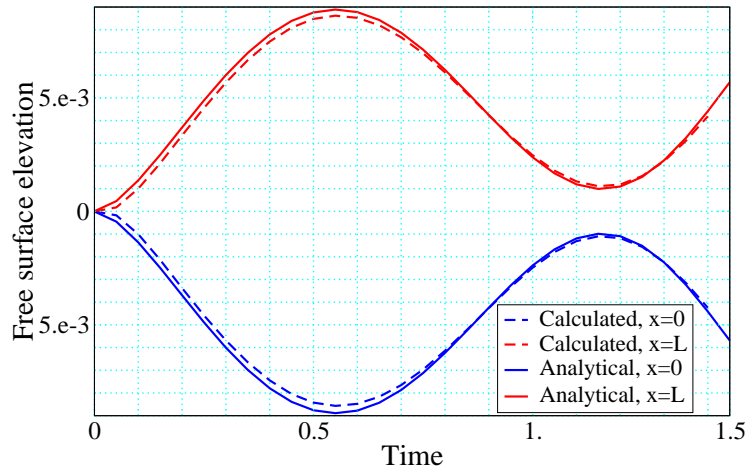}} \end{center}
\caption{Sloshing in cavity: analytical solution and numerical solution with $\mu= \rho/1000$.
\label{ball_num_rho1000}}
\end{figure}

\begin{figure}[htb]
\begin{center} \scalebox{0.6}{\includegraphics*{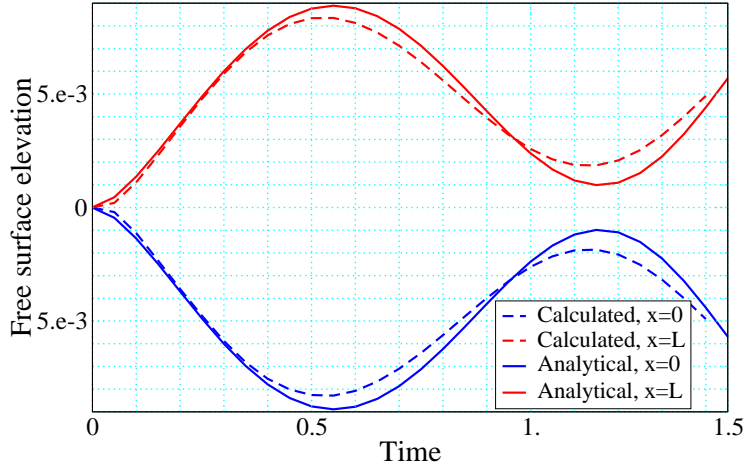}} \end{center}
\caption{Sloshing in cavity: analytical solution and numerical solution with $\mu= \rho/10000$.
\label{ball_num_rho10000}}
\end{figure}

\medskip
The numerical results are reported on figure \ref{ball_num_rho100} ($\mu= \rho/100$), figure \ref{ball_num_rho1000} ($\mu= \rho/1000$), and figure \ref{ball_num_rho10000} ($\mu= \rho/10000$) respectively. 
Comparing the obtained shape for the interface with the analytical solution, we observe that the numerical solution is closer to the analytical one with $\mu= \rho/1000$ than with $\mu= \rho/100$, certainly because the fluid is too viscous in this latter case.
More surprisingly, when reducing the viscosity to $\mu= \rho/10000$, the numerical solution also becomes less accurate.
Our explanation is that, to obtain a good solution, it is necessary to respect a balance between approaching the physical problem (which, in this case, would suggest $\mu=0$) and keeping sufficient coercivity to ensure a reasonable convergence of the numerical approximation (which, on the contrary, requires a high value for the viscosity).
With a more refined mesh, viscosity thus probably could be decreased, and the solution be closer to the analytical one.
However, with this mesh already, results seem to be rather more accurate as those available in the litterature \cite{cha-04-mod}.

\subsection{Bubble column}

\begin{figure}
\begin{center}
\psfrag*{L}{{\Huge $L$}}
\psfrag*{H}{{\Huge $H$}}
\psfrag*{h}{{\Huge $h$}}
\psfrag*{q}{{\Huge $q$}}
\begin{tabular}[t]{cccc}
\scalebox{0.33}{\includegraphics*[0.cm,0.cm][15cm,25cm]{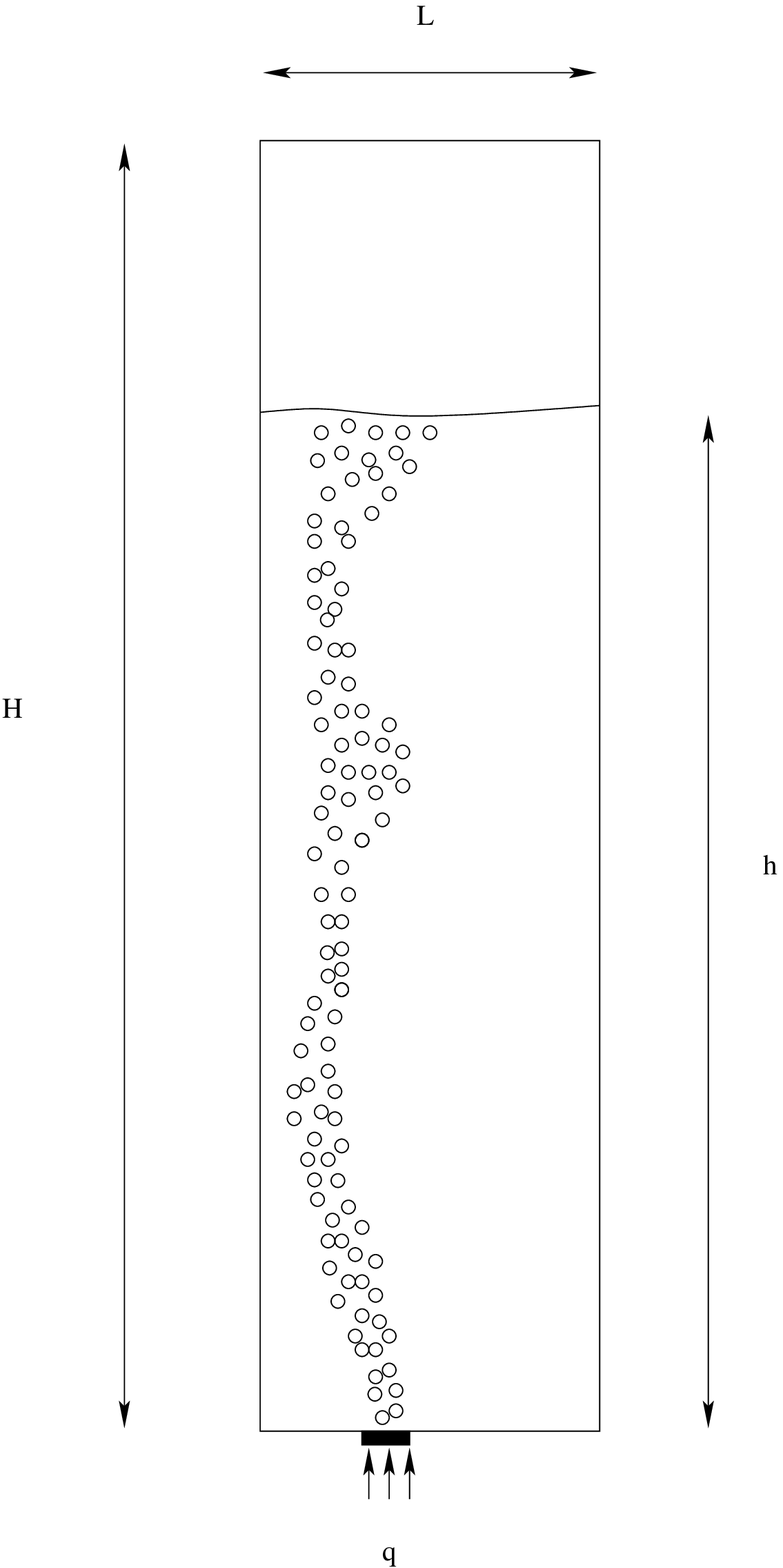}}
&
\scalebox{0.35}{\includegraphics*[8.cm,3.cm][15cm,25cm]{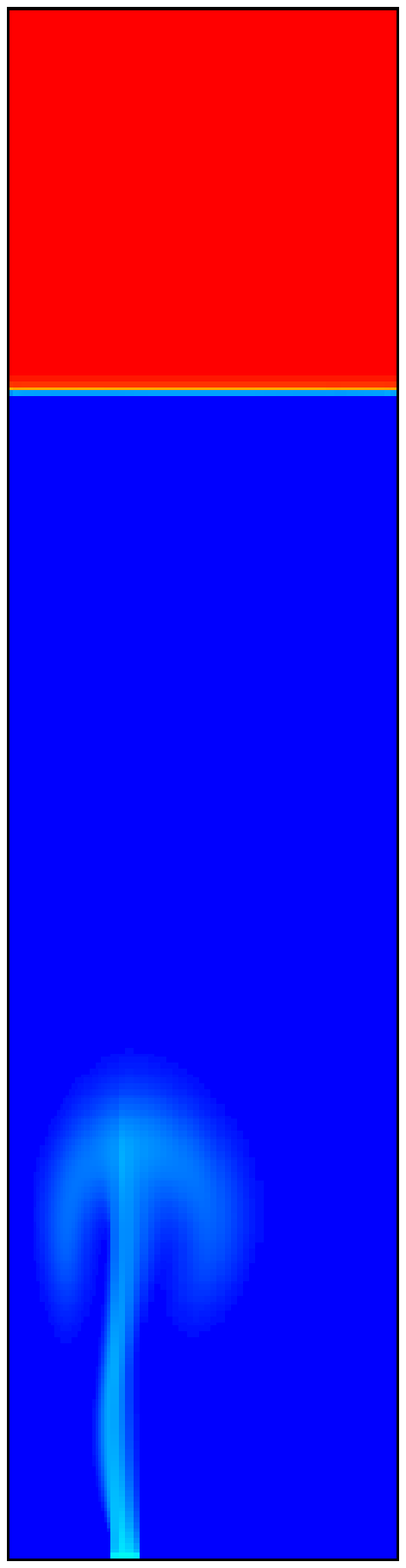}}
&
\scalebox{0.35}{\includegraphics*[8.cm,3.cm][15cm,25cm]{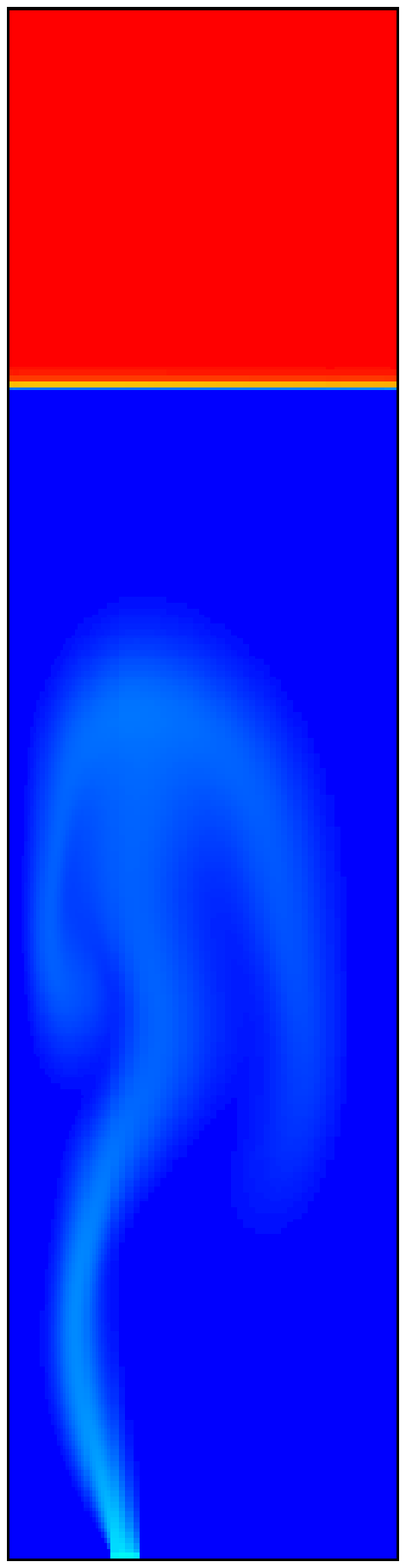}}
&
\scalebox{0.35}{\includegraphics*[8.cm,3.cm][22cm,27cm]{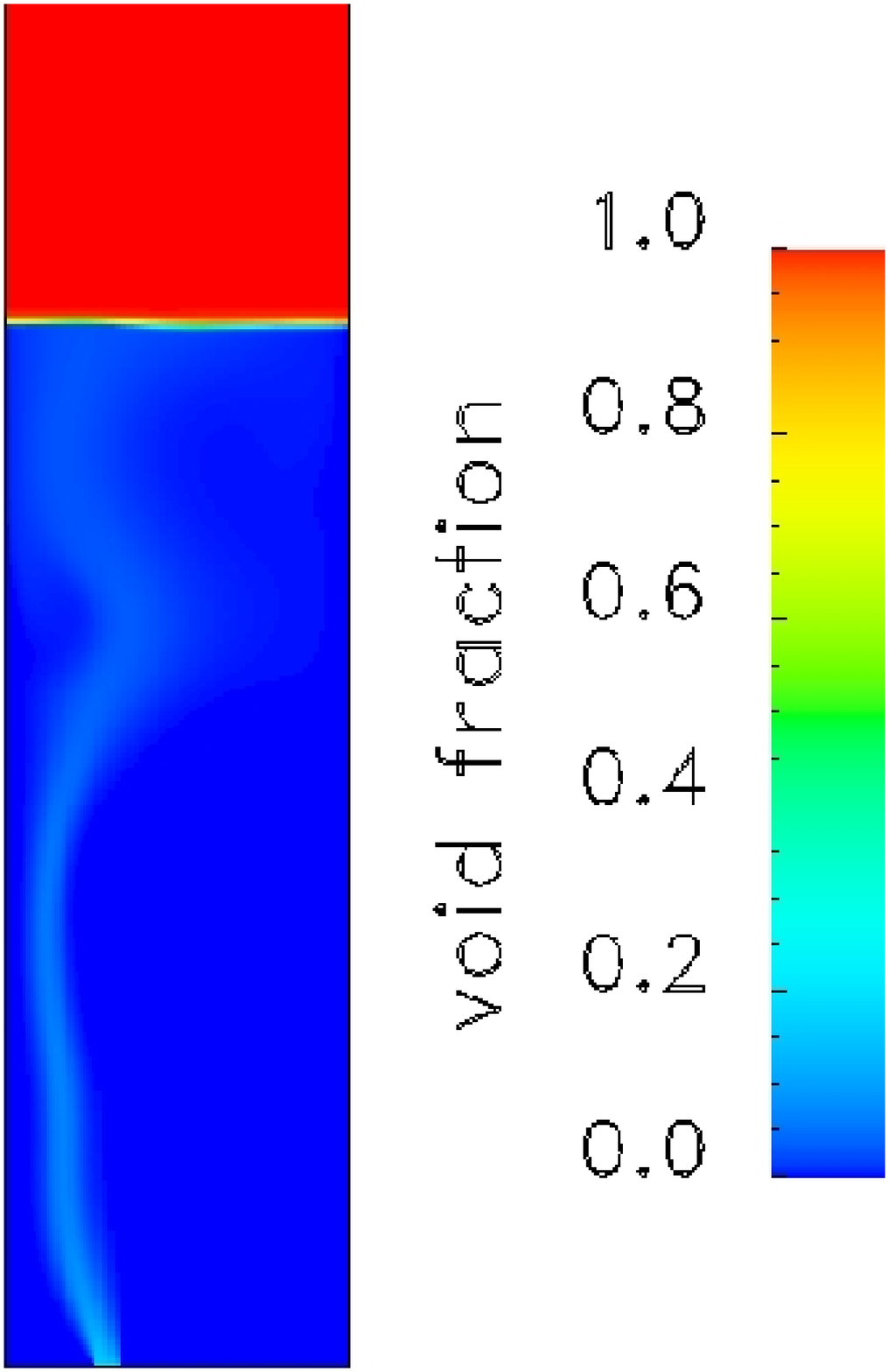}}
\end{tabular}
\end{center}
\caption{Bubble column: geometry of the problem and void fraction at times $2s$, $4s$ and $40s$.
\label{bubble_res}}
\end{figure}

We address in this section a classical benchmark for diphasic flow solvers, namely the flow in a pseudo two dimensional bubble column  investigated experimentally by  Becker \etal \cite{bec-94-gas}.
The apparatus has a rectangular cross section with the following dimensions: its width is $L=50\,cm$, its depth is $8\,cm$ and it is $H=200\,cm$ high (see figure \ref{bubble_res}).
It is filled with water up to the height $h=150\,cm$.
A gas sparger, positioned $15\,cm$ from the left wall, is used to introduce an air flow of $q=8\,l.min$ into the system. The circular sparger has a diameter of $40\,mm$ and a pore size of $40\,\mu m$.
Several liquid circulation cells can be observed in the column, the location and size of which continously change.
The bubble swarm is influenced by these vortices and therefore rises in a meander-like way.
The direction of its lower part is stable and directed towards the nearest sidewall; its upper part changes its shape and location in a quasiperiodic way, according to transient liquid circulations \cite{sok-99-app}.

\medskip
To simulate this experiment, we choose the following data.
The boundary conditions are defined at the inlet as follows:
\[
u_{imp}=\frac{q}{S\,\alpha_{g,{\rm imp}}}
\]
where $S$ is the gas inlet area and $\alpha_{g,{\rm imp}}=1$ is the void fraction imposed at the inlet.
Along the walls and at the outlet of the column,  homogeneous Dirichlet conditions are used for the velocity. 
Initial conditions are set to $u=0 \,m.s^{-1}$ and $p=p_0$ where $p_0=10^5 \, Pa$ is the ambiant pressure.
The density of the liquid is $\rho_{\ell}=1000 \, kg.m^{-3}$; the gas obeys an ideal gas equation of state $\rho_g = p/a^2$, where $a^2$ is such that $\rho_g = 1.2 \, kg.m^{-3}$ at $p=10^5 \, Pa$. 
The diffusion coefficient $D$ is set to zero, the drift velocity is constant and given by $u_r=(0,\,0.2)^t \,m.s^{-1}$.

\medskip
For this test case, we use a regular meshing composed of rectangular cells (with the Rannacher-Turek element) with 76 meshes in the horizontal direction, from which 4 for the gas inlet, and 300 in the vertical one.
Calculations with time steps up to $\dt=10^{-1}\, s$ have been performed, observing that smaller time steps yield a thinner free surface.

\medskip
The viscosity is a parameter difficult to adjust, since, in this simulation which is based on the system of equations governing a laminar flow, it must represent in some way the turbulent diffusion, \ie\ the effects of fluctuations of the flow at microscopic scales, which may originate from the usual turbulence phenomena (sometimes termed "monophasic turbulence") and from the perturbation of the velocity field due to the motion of the bubbles (sometimes termed "diphasic turbulence").
Calculations with a viscosity ranging from $\mu=10^{-3}\,Pa.s$ to $\mu=10^{2}\,Pa.s$ have been performed.
With smaller viscosities, we observe more oscillations of the free surface, the bubble swarm reaches the free surface faster and is farther from the sidewall.
 
\medskip
Finally, the numerical results obtained with $\dt = 10^{-2}\, s$ and a viscosity of $\mu=1\,Pa.s$ are reported on figure \ref{bubble_res}.
With this value of the viscosity and these mesh and time steps, numerical convergence seems to be reached, at least visually.
One can observe the stability and the thinness of the free surface.
Results qualitatively reproduce the expected behaviour, which is the best we can hope with the rather crude modelling of turbulence which we adopted.

\section{Conclusion}

In this paper, we adress the drift-flux model, which, for isothermal flows, consists in a system of three balance equations, namely the overall mass, the gas mass and the momentum balance complemented by an equation of state and a phenomenologic relation for the drift velocity.

\medskip
For this problem, we develop a pressure correction scheme combining finite element and finite volume discretizations, which enjoys the following properties.
First, the existence of a solution to each step of the algorithm is proven.
Then, essential stability features of the continuous problem still hold at the discrete level: the unknowns are kept within their physical bounds (in particular, the gas mass fraction remains in the $[0,1]$ interval); in the homogeneous case (\ie\ when the drift velocity vanishes), the discrete entropy of the system decreases; in addition, when using for the drift velocity the Darcy-like relation suggested in \cite{gui-07-ada}, the drift term becomes dissipative.
Since, when the density is constant, this fractional step algorithm degenerates to an usual incremental projection method based on an {\em inf-sup} stable approximation, stability can be expected in the zero Mach number limit.
Finally, the present algorithm preserves a constant pressure and a constant velocity through moving interfaces between phases (\ie\ contact discontinuities of the underlying hyperbolic system).
To achieve this latter goal, the key ingredient is to couple the mass balance and the transport terms of the gas mass balance in an original pressure correction step.

\medskip
We chose in this paper to only consider the case of a constant density liquid phase and of a gaseous phase obeying the ideal gas law.
Dealing with a more general barotropic gas phase is certainly the simplest generalization, but the present theory also seems to extend to the case of a compressible fluid with minor modifications: for the stability study, essentially, the expression for the volumetric free energy of the mixture should be replaced by the usual expression applying when both phases are compressible, see for instance \cite{gui-07-ada}; the existence theory would probably be simpler, since an upper bound for the density would provide in this case an estimate for the pressure.
Returning to the case of an incompressible fluid, extending the present theory to deal with pure liquid zones appears on the contrary to be a difficult task, since the role played by the pressure in such a system seems to deserve some clarifications.

\medskip
Numerical tests show a near-first-order convergence in space and time, consistent with the implemented discretization: first order backward Euler method in time and standard upwinding of the convection terms in the mass and gas mass fraction balance equations.
With respect to this latter point, using more accurate space discretization (typically, MUSCL-like techniques) should certainly be desirable.

\medskip
To assess the robustness of this algorithm, various numerical tests have been performed.
They show in particular that free surface flows are computed without any instability, keeping a rather sharp interface throughout the computation.
In addition, pure monophasic liquid zones are supported, although, as already mentioned, this case remains beyond the scope of the theory developped here.
This scheme is now implemented in the ISIS code developped at IRSN and daily used for industrial applications.


\appendix

\section{Existence of a solution to a class of discrete diphasic problems}\label{sec:existence_z}

We address in this section the following abstract discrete problem:
\begin{equation}
\left| \begin{array}{l} \displaystyle
a(u, \varphi_{\edge}^{(i)})
- \int_{\Omega,h} p\ \nabla \cdot \varphi_{\edge}^{(i)}{\rm d}x = \int_\Omega f \cdot \varphi_{\edge}^{(i)} {\rm d}x,
\hfill
\forall \edge\in\edgesint,\ 1 \leq i \leq d
\\[2ex] \displaystyle
\frac{|K|}{\dt} \left[ \varrho^{\,p,z}(p_K,\, z_K)-\varrho^{\,p,z}(p_K^\ast,\, z^\ast_K) \right]
+ \sum_{\edge=K|L} \vsp \varrho^{\,p,z}(p_K,z_K) - \vsm \varrho^{\,p,z}(p_L,z_L)=0,
\qquad \hfill
\forall K \in \mesh
\\[2ex] \displaystyle
\frac{|K|}{\dt} (z_K-z^\ast_K) + \sum_{\edge=K|L} \vsp\, z_K - \vsm\ z_L=0,
\hfill
\forall K \in \mesh
\end{array}\right .
\label{pbdiph_disc}\end{equation}
This problem is supposed to be obtained from (part of) a continuous problem by a space discretization combining Rannacher-Turek or Crouzeix-Raviart finite elements and finite volumes; notations related to discrete quantities are given in section \ref{sec:spacedisc} and are not recalled here.  
The bilinear form $a(\cdot,\cdot)$ is just assumed to be such that $\norma{u}=[a(u,u)]^{1/2}$ defines a norm over the discrete space $W_h$.
The quantities $(\vsp)^{n+1}$ and $(\vsm)^{n+1}$ stands respectively for $\max (\vs^{n+1},\ 0)$ and $-\min (\vs^{n+1},\ 0)$ with $\vs^{n+1}=|\edge|\, u_\edge^{n+1} \cdot n_{KL}$.
Note that, in the last two equations, the flux summation excludes the external edges, which implicitely expresses the fact that the velocity is supposed to vanish on the boundary.

\medskip
This system must be completed by three equations of state.
The first two ones giving the liquid density $\rho_\ell$ and the gas density $\rho_g$ as a function of the pressure: we suppose here that the density of the liquid is constant and that the gas obeys the equation of state of ideal gases, which, for the sake of conciseness, we suppose here to be simply $\rho_g=p$.
The last equation relates the mixture density $\rho$ with the gas mass fraction $y$ or the gas partial density $z=\rho\,y$ and the phases density, and may take the three following forms:
\begin{equation}
p=\wp(\rho,z)=\frac{z \ \rho_\ell}{z+\rho_\ell-\rho} \qquad \qquad
\rho=\varrho^{\,p,z}(p,z)= \rho_\ell + z\ (1 - \frac {\rho_\ell} p) \qquad \qquad
\rho=\varrho^{\,p,y}(p,y)=\frac 1 {\displaystyle \frac{1-y}{\rho_\ell} + \frac y p}
\label{rhomixture}\end{equation}
These three relations are equivalent as soon as the following assumptions for the unknowns of this system are satisfied:
\begin{equation}
\rho>0, \qquad
p>0, \qquad 
z > 0 \qquad \mbox{and} \quad
0<y\leq 1
\label{phys_bounds}\end{equation}  
These assumptions are natural, excepted the hypothesis that $y$ or $z$ does not vanish, which excludes the existence of purely liquid zones.
This latter assumption is assumed to hold for the initial quantities, \ie\ we suppose that:
\begin{equation}
\forall K \in \mesh, \qquad y_K^\ast=\frac{z_K^\ast}{\rho_K^\ast} \in (0,1]
\label{ystar}\end{equation}
where $\rho_K^\ast=\varrho^{\,p,z}(p_K^\ast,\, z^\ast_K)$.

\medskip
Our aim in this section is to prove that there exists a solution to system \eqref{pbdiph_disc} complemented with one of the relations of \eqref{rhomixture}, under the assumption \eqref{ystar}, and that any such solution satisfies the inequalities \eqref{phys_bounds}.

\bigskip
We begin this section by two preliminary lemmas.

\begin{lmm}
Let $(x_K^\ast)_{K\in \mesh}$ and $(x_K)_{K\in \mesh}$ be two families of real numbers satisfying the following set of equations:
\[
\forall K \in \mesh, \qquad
\frac{|K|}{\dt}\,(x_K-x^\ast_K) + \sum_{\edge=K|L}\,\vsp\,x_K-\vsm\,x_L=0
\]
We suppose that, $\forall K \in \mesh,\ x_K^\ast >0$.
Let $\norminf{\nabla_h \cdot u}$ be defined by:
\[
\norminf{\nabla_h \cdot u}=\max_{K\in\mesh} \left[ 0,\frac{1}{|K|} \sum_{\edge=K|L} \vs \right]
\]
Then, $\forall K\in \mesh$, $x_K$ satisfies:
\[
\frac{\displaystyle \min_{K\in\mesh} x_K^\ast}{1 + \dt \norminf{\nabla_h \cdot u}} 
\leq x_K \leq 
\frac{1}{\displaystyle \min_{K\in\mesh}|K|}\sum_{K\in \mesh} |K|\ x^\ast_K
\]
\label{bound_rho_z}\end{lmm}

\begin{proof}
The first inequality follows from an application of the discrete maximum principle lemma which can be found in \cite{gal-07-an} (lemma 2.5, section 2.3).
The second one then follows from the fact that, by conservativity, $\sum_{K\in\mesh} x_K = \sum_{K\in\mesh} x_K^\ast$, remarking that, by the preceding relation, the values $x_K$, for $K\in\mesh$, are all positive.
\end{proof}

The proof of the following result can be found in \cite{lar-91-how}.
\begin{lmm}
Let $(\rho_K^\ast)_{K\in \mesh}$, $(x_K^\ast)_{K\in \mesh}$, $(\rho_K)_{K\in \mesh}$ and $(x_K)_{K\in \mesh}$ be four families of real numbers satisfying the following set of equations:
\[
\forall K \in \mesh, \qquad
\frac{|K|}{\dt}\, (\rho_K\, x_K-\rho_K^\ast\, x^\ast_K) + \sum_{\edge=K|L}\,\vsp\,\rho_K\,x_K-\vsm\,\rho_L\,x_L=0
\]
We suppose that, $\forall K \in \mesh,\ \rho_K^\ast>0$, $\rho_K >0$ and:
\[
\forall K \in \mesh, \qquad
\frac{|K|}{\dt}\,(\rho_K-\rho_K^\ast) + \sum_{\edge=K|L}\,\vsp\,\rho_K-\vsm\,\rho_L=0
\]
Then the following discrete maximum principle holds:
\[
\forall K \in \mesh,\qquad \min_{L\in \mesh} x_L^\ast \leq x_K \leq \max_{L\in \mesh} x_L^\ast
\]
\label{larrout'}\end{lmm}

We now state the abstract theorem which will be used hereafter; this result follows from standard arguments of the topological degree theory (see \cite{deimling} for an overview of the theory and \eg\ \cite{eym-98-an, gal-07-an} for other uses in the same objective as here, namely the proof of existence of a solution to a numerical scheme).
%
%
\begin{thrm}[A result from the topological degree theory]\label{degre} 
Let $N$ and $M$ be two positive integers and $V$ be defined as follows:
\[
V=\{ (x,y,z) \in \xR^N \times \xR^M \times \xR^M \mbox{ such that } y>0 \mbox{ and } z>0\}
\]
where, for any real number $c$ and vector $y$, the notation $y>c$ means that each component of $y$ is greater than $c$.
Let $b\in \xR^N \times \xR^M \times \xR^M$ and $f(\cdot)$ and $F(\cdot,\cdot)$ be two continuous functions respectively from $V$ and $V\times[0,1]$ to $\xR^N \times \xR^M \times \xR^M$ satisfying:
\begin{enumerate}
\item[(i)]$F(\cdot,\,1)=f(\cdot)$;
\item[(ii)] $\forall \theta \in [0,1]$, if an element $v$ of $\bar {\cal O}$ (the closure of ${\cal O}$) is such that $F(v,\theta)=b$, then $v \in {\cal O}$, where ${\cal O}$ is defined as follows:
\[
{\cal O}=\{ (x,y,z) \in \xR^N \times \xR^M \times \xR^M \mbox{ s.t. } \Vert x \Vert < M \mbox{ and } \epsilon< y < M
\mbox{ and } \epsilon < z < M\}
\]
with $M$ and $\epsilon$ two positive constants and $\Vert \cdot \Vert$ a norm defined over $\xR^N$;
\item[(iii)] the topological degree of $F(\cdot,0)$ with respect to $b$ and ${\cal O}$ is equal to $d_0 \neq 0$.
\end{enumerate}
Then the topological degree of $F(\cdot,1)$ with respect to $b$ and ${\cal O}$ is also equal to $d_0 \neq 0$; consequently, there exists at least a solution $v\in {\cal O}$ such that $f(v)=b$.
\end{thrm}

\smallskip
We are now in position to prove the existence of a solution to the considered discrete system.
%
%
\begin{thrm}[Existence of a solution]
Under the assumption \eqref{ystar}, the nonlinear system \eqref{pbdiph_disc} complemented with the relation \eqref{rhomixture} admits at least one solution, and any possible solution is such that:
\[
\forall K \in \mesh,\qquad \rho_K >0, \qquad z_K >0, \qquad 0 < y_K=\frac{z_K}{\rho_K} \leq 1, \qquad p_K >0
\]
\label{existence_pxz}\end{thrm}

\begin{proof}
This proof makes use of theorem \ref{degre} twice, by linking the initial problem to a linear one through two successive homotopies.
Let $N=d\ {\rm card}(\edgesint)$ and $M={\rm card}(\mesh)$; we identify the finite element space of discrete velocity with $\xR^N$ and the finite volume space of pressure and partial density with $\xR^M$.
Let $V$ be defined by $V=\{ (u,p,z) \in \xR^N\times \xR^M\times \xR^M \mbox{ such that } p>0 \mbox{ and } z>0\}$.

\vspace{3ex}
\noindent \underline{Step 1: first homotopy}
\\[1ex]
We consider the function $F:\ V\times [0,1] \rightarrow \xR^N\times \xR^M\times \xR^M$ given by:
\begin{equation}
\begin{array}{ll}
F(u,p,z,\theta)=
\left| \begin{array}{l} \displaystyle
v_{\edge,i}=
a(u,\varphi_\edge^{(i)})
- \int_{\Omega,h} p\ \nabla \cdot \varphi_\edge^{(i)} \, {\rm d}x- \int_\Omega f \cdot \varphi_\edge^{(i)} \, {\rm d}x,
\hfill
\edge \in \edgesint,\ 1 \leq i \leq d
\\[2ex] \displaystyle
q_K=\frac{|K|}{\dt} \left[\varrho^{\,p,z}_\theta(p_K,z_K)-\varrho^{\,p,z}_\theta(p_K^\ast,z_K^\ast)\right]
\\[2ex] \displaystyle \hspace{15ex}
+ \sum_{\edge=K|L} \vsp\, \varrho^{\,p,z}_\theta(p_K,z_K) - \vsm\, \varrho^{\,p,z}_\theta(p_L,z_L),
\hspace{10ex}
K \in \mesh
\\ [2ex] \displaystyle
s_K=\frac{|K|}{\dt} \left[ z_K- \varrho^{\,p,z}_\theta(p_K^\ast,z_K^\ast) \, y_K^\ast \right] 
+ \sum_{\edge=K|L} \vsp\, z_K - \vsm\, z_L,
\hfill
K \in \mesh
\end{array} \right .
\end{array}
\label{eqalpha} \end{equation}
where the function $\varrho^{\,p,z}_\theta(\cdot,\cdot)$ is implicitely defined by the following relation:
\[
\varrho^{\,p,z}_\theta (p,z)=\varrho^{\,p,y}_\theta (p,y)=\frac 1 {\displaystyle \frac{1-y}{\varrho_{\ell,\theta}(p)} + \frac y p}
\qquad \mbox{with} \quad
\varrho_{\ell,\theta}(p)=\frac 1 {\displaystyle \frac \theta {\rho_\ell} + \frac {1-\theta} p}
\quad \mbox{ and } z=\rho\,y
\]
Note that this definition makes sense (\ie\ using $z=\rho y$, the function $\varrho^{\,p,z}(\cdot,\cdot)$ can be explicitely computed from the expression of $\varrho^{\,p,y}(\cdot,\cdot)$) as soon as $p>0$, and thus for any $(u,p,z) \in V$.

\medskip
Problem $F(u,p,z,1)=0$ is exactly the same as system \eqref{pbdiph_disc}.

\medskip
Let $\epsilon$ and $M$ be two positive real numbers, and $\cal O$ be defined by:
\[
{\cal O}=\{ (u,p,z) \in \xR^N \times \xR^M \times \xR^M \mbox{ s.t. } \norma{u} < M,\ \epsilon< p < M
\mbox{ and } \epsilon < z < M\}
\]
We now suppose that $(u,p,z) \in \bar {\cal O}$ (and thus, in particular, $p \geq \epsilon$) and that $F(u,p,z,\theta)=0$ and provide estimates for $(u,p,z)$. 

\medskip
We begin by the following elementary bound, which is useful throughout the proof.
From the definition of $\varrho_{\ell,\theta}(p)$, we observe that $\min(\rho_\ell,p) \leq \varrho_{\ell,\theta}(p) \leq \max(\rho_\ell,p)$.
By the same way, provided that $y \in [0,1]$, $\min(\varrho_{\ell,\theta}(p),p) \leq \varrho^{\,p,z}_\theta(p,z) \leq \max(\varrho_{\ell,\theta}(p),p)$.
Hence, $\min(\rho_\ell,p) \leq \varrho^{\,p,z}_\theta(p,z) \leq \max(\rho_\ell,p)$ and, thanks to assumption \eqref{ystar}:
\[
\forall \theta \in [0,1],\ \forall K \in \mesh, \qquad \varrho^{\,p,z}_\theta(p_K^\ast,z_K^\ast)\leq \bar \rho^\ast
\qquad \mbox{with} \quad \bar \rho^\ast=\max \left[ (\max_{K\in\mesh} p_K^\ast),\ \rho_\ell \right]
\] 

\bigskip
\noindent \underline{\textit{Step 1.1}}: $\norma{\cdot}$ estimate for the velocity.
\\[1ex]
Let us first recast the equation of state of the mixture under a more convenient form.
Substituting its definition for $\varrho_{\ell,\theta}(p)$ in $\varrho^{\,p,y}_\theta(p,y)$, we get:
\begin{equation}
\rho=\frac 1 {\displaystyle \frac{\theta \, (1-y)}{\varrho_\ell} + \frac {y+(1-\theta)(1-y)} p}
= \frac 1 {\displaystyle \frac{1-y'}{\varrho_\ell} + \frac {y'} p}
\label{rho_mixt_alpha}\end{equation}
with $y'(y,\theta)=y+(1-\theta)(1-y)$.
Then, taking $y=z/\varrho^{\,p,z}_\theta(p,z)$ as unknown in the third equation of $F(u,p,z,\theta)=0$, we get, for any $K \in \mesh$:
\[
\frac{|K|} \dt \, \left[ \varrho^{\,p,z}_\theta(p_K,z_K)\, y_K- \varrho^{\,p,z}_\theta(p_K^\ast,z_K^\ast)\, y_K^\ast \right]
 + \sum_{\edge=K|L} \vsp\, \varrho^{\,p,z}_\theta(p_K,z_K)\, y_K - \vsm\, \varrho^{\,p,z}_\theta(p_L,z_L)\, y_L=0
\]
As, by the second equation of $F(u,p,z,\theta)=0$, this relation vanishes for the constant function $y_K=1,\ \forall K \in \mesh$, we also obtain:
\begin{equation}
\frac{|K|} \dt \, \left[\varrho^{\,p,z}_\theta(p_K,z_K)\, y'_K- \varrho^{\,p,z}_\theta(p_K^\ast,z_K^\ast)\, (y')_K^\ast \right]
 + \sum_{\edge=K|L} \vsp\, \varrho^{\,p,z}_\theta(p_K,z_K)\, y'_K - \vsm\, \varrho^{\,p,z}_\theta(p_L,z_L)\, y'_L=0
\label{pbz'}\end{equation}
where, by assumption \eqref{ystar}, $\forall K \in \mesh,\ (y')_K^\ast=y_K^\ast+(1-\theta)(1-y_K^\ast) \in (1-\theta + \theta \underbar y^\ast,1]$, with $\underbar y^\ast=\min_{K\in\mesh} y_K^\ast$.
We thus obtain a new problem, which keeps the structure of system \eqref{pbdiph_disc}, with the same equation of state (\ie\ relation \eqref{rho_mixt_alpha}) and just a modified initial value for $z$ (\ie\ $z_K^\ast$ changed to $\varrho^{\,p,z}_\theta(p_K^\ast,z_K^\ast)\, (y')_K^\ast$).
The unknown $p$ is still an unknown of this new problem, and we thus have by assumption $p\geq 0$.
In addition, by lemma \ref{larrout'}, any solution of this new problem is such that the gas mass fraction verifies $1-\theta + \theta \underbar y^\ast < y \leq 1$, and thus the density and the gas partial density are positive.
The unknowns thus belong to the domain where the free energy is correctly defined, and theorem \ref{pot_el_pxz} applies.
Multiplying the first equation of $F(u,p,z,\theta)=0$ by $u_{\edge,i}$, summing over $\edge \in \edgesint,\ 1 \leq i \leq d$ and using Young's inequality thus yields:
\[
\frac 1 2 \normad{u} + \frac 1 \dt \sum_{K\in\mesh} |K|\ \varrho^{\,p,z}_\theta(p_K,z_K)\, y'_K\ \log(p_K) \leq 
\frac 1 \dt \sum_{K\in\mesh} |K|\ \varrho^{\,p,z}_\theta(p^\ast_K,z^\ast_K)\, (y')^\ast_K\ \log(p^\ast_K) + \frac 1 2 \normmad{f}
\]
where $\normma{\cdot}$ stands for the dual norm of $\norma{\cdot}$ with respect to the $\xLtwo$ inner product.
The summation at the right hand side of this relation is bounded by $(1/\dt)\ |\Omega|\ \bar \rho^\ast \log(\bar p^\ast)$ where $\bar p^\ast=\max_{K\in\mesh}p_K^\ast$.
By conservativity of equation \eqref{pbz'}, $\sum_{K\in\mesh} |K|\ \varrho_\theta(p^\ast_K,z^\ast_K)\, (y')^\ast_K \leq |\Omega|\ \bar \rho^\ast$.
Since, by assumption, $p \geq \epsilon$, we thus get:
\[
\normad{u} \leq \frac 2 \dt\ |\Omega|\ \bar \rho^\ast\ |\log (\epsilon)| + \frac 2 \dt \ |\Omega|\ \bar \rho^\ast \log(\bar p^\ast) + \normmad{f}
\]
For $\epsilon>0$ small enough, we thus have:
\begin{equation}
\norma{u} \leq c_1 \ |\log (\epsilon)|^{1/2}
\label{est_u}\end{equation}
where, in this relation and throughout the proof, we denote by $c_i$ a real number only depending on the data of the problem, \ie\ $\Omega$, $\bar \rho^\ast$, $\bar p^\ast$, $f$, $a(\cdot,\cdot)$, $\dt$ and the mesh, and the expression "$\epsilon$ small enough" stands for $\epsilon < c'_1$ where $c'_1$ is a positive real number itself only depending on the data.

\bigskip
\noindent \underline{\textit{Step 1.2}}: $\xLinfty$ estimates for $z$.
\\[1ex]
By equivalence of the norms over finite dimensional spaces, inequality \eqref{est_u} also yields a bound for $u$ in the $\xLinfty$ norm and, finally, for $\norminf{\nabla_h \cdot u}$:
\[
\norminf{\nabla_h \cdot u} \leq c_2 \ |\log (\epsilon)|^{1/2}
\]
By lemma \ref{bound_rho_z}, we thus get from the third relation of the system $F(u,p,z,\theta)=0$, still for $\epsilon$ small enough:
\begin{equation}
z \geq c_3 \ |\log (\epsilon)|^{-1/2}
\label{est_z_inf}\end{equation}
On the other hand, we get from the same relation by conservativity:
\begin{equation}
z \leq c_4
\label{est_z_sup}\end{equation}

\bigskip
\noindent \underline{\textit{Step 1.3}}: $\xLinfty$ estimates for $p$.
\\[1ex]
From the first relation of \eqref{rhomixture}, using the bounds for $z$, we get:
\begin{equation}
p \geq c_5 \ |\log (\epsilon)|^{-1/2}
\label{est_p_inf}\end{equation}
To obtain an upper bound for $p$, we first remark that, as the considered spatial discretization satifies a discrete {\em inf-sup} condition, a bound for $u$ provides a bound for $p - m(p)$ where $m(p)$ stands for the mean value of $p$.
By equivalence of norms on finite dimensional spaces, we can choose to express this bound in the seminorm defined by $\forall q \in L_h,\ |q|_{1,1,h}=\sum_{\edge \in \edgesint\ (\edge=K|L)} |q_K-q_L|$.
With this semi-norm, the mean value of $p$ disappears, and we get for $\epsilon$ small enough:
\begin{equation}
\sum_{\edge \in \edgesint\ (\edge=K|L)} |p_K-p_L| \leq c_6\ |\log (\epsilon)|^{1/2}
\label{est_p_sup_prel}\end{equation}
An upper bound for $p$ in one cell of the mesh, say $K_0$, would then provide an upper bound for $p$, since, for any $K \in \mesh$, it is possible to build a path from $K_0$ to $K$ crossing each internal edge at most once.
To obtain such an estimate, we follow the following idea.
If the pressure is somewhere lower than $\rho_\ell$, we are done; otherwise, with the chosen equation of state $\varrho^{\,p,z}_\theta(\cdot,\cdot)$, when $\theta$ varies, the liquid is everywhere denser than for $\theta=1$ and we are going to show that, even if its total mass also increases, the volume that it occupies is lower than for $\theta=1$.
Hence, the remaining volume for the gas is bounded away from zero, and, by conservation of the gas mass, the pressure cannot blow up everywhere.
First, we need to introduce the phase volumetric fractions.
The equation of state \eqref{rhomixture} can be written as:
\[
\frac \rho p + \frac{\rho -z}{\rho_\ell}=1
\]
and, as $\rho -z = \rho\ (1-y)$ and $y \leq 1$, both fractions at the left hand side of this relation are non-negative.
We may thus define $\alpha_g \in [0,1]$ and $\alpha_\ell \in [0,1]$, referred to as the gas and liquid volume fraction respectively, by:
\[
\alpha_g= \frac \rho p \qquad \qquad \alpha_\ell=\frac{\rho -z}{\rho_\ell}
\]
Note that $\alpha_g + \alpha_\ell=1$.
Combining the second and the third relation of the system $F(u,p,z,\theta)=0$, summing over the control volumes of the mesh and remarking that the fluxes cancel by conservativity, we get:
\begin{equation}
\sum_{K \in \mesh} |K| \ (\alpha_\ell)_K = 
\sum_{K \in \mesh} |K| \ \frac{(1-y_K^\ast)\ \varrho^{\,p,z}_\theta(p_K^\ast,z_K^\ast)}{\rho_\ell}
\label{alpha1}\end{equation}
Let us denote by $(\alpha^\ast_\ell)_{K,1}$ the liquid void fraction with the equation of state corresponding to $\theta=1$:
\[
(\alpha^\ast_\ell)_{K,1}=\frac {1 - y_K^\ast} {\displaystyle \rho_\ell \left[ \frac{y_K^\ast}{p_K^\ast}+ \frac{1-y_K^\ast}{\rho_\ell}\right]}
\]
Exploiting the form \eqref{rho_mixt_alpha} of the equation of state for $\theta \neq 0$, we obtain from relation \eqref{alpha1}:
\[
\sum_{K \in \mesh} |K| \ (\alpha_\ell)_K = \sum_{K \in \mesh} |K| (\alpha^\ast_\ell)_{K,1}
\ \frac{\displaystyle \frac{y_K^\ast}{p_K^\ast}+ \frac{1-y_K^\ast}{\rho_\ell}}
{\displaystyle \frac{(y')_K^\ast}{p_K^\ast}+ \frac{1-(y')_K^\ast}{\rho_\ell}}
\qquad \mbox{ with} \quad
(y')_K^\ast=y_K^\ast + \theta (1-y_K^\ast)
\]
If we suppose that $p_K \geq \rho_\ell$, the fraction in the above equation is bounded by $1$: indeed, both the numerator and the denominator are harmonic averages of $p_K^\ast$ and $\rho_\ell$, the weight associated to $p_K^\ast$ being larger in the denominator, since $(y')_K^\ast$ is closer to $1$ than $y_K^\ast$.
We thus get:
\[
\sum_{K \in \mesh} |K| \ (\alpha_g)_K \geq c_7 = |\Omega| - \sum_{K \in \mesh} |K| (\alpha_\ell)_{K,1}
\]
where $c_7$ is positive by assumption, since $\forall K \in\mesh,\ y_K^\ast > 0$.
Thus there exists $K_0 \in \mesh$ such that $(\alpha_g)_{K_0} \geq c_8=c_7/|\Omega|$.
On the other hand, we have, still by conservativity:
\[
\sum_{K \in \mesh} |K| \ (\alpha_g)_K \ p_K= \sum_{K \in \mesh} |K|\ z_K = \sum_{K \in \mesh} |K|\ z^\ast_K
= \sum_{K \in \mesh} |K|\ \varrho^{\,p,z}_\theta(p_K^\ast,z_K^\ast) \, y^\ast_K
\leq |\Omega|\ \bar \rho^\ast
\]
We thus get, since all the $(\alpha_g)_K$ and $p_k$ are non-negative:
\[
(\alpha_g)_{K_0} \ p_{K_0} \leq \frac{|\Omega|}{|K_0|}\ \bar \rho^\ast
\]
and thus, as $(\alpha_g)_{K_0}$ is bounded by below, the pressure is bounded by a quantity only depending on the data.
As a consequence, for $\epsilon$ small enough:
\begin{equation}
p \leq c_9\ |\log (\epsilon)|^{1/2}
\label{est_p_sup}\end{equation}

\clearpage
\noindent \underline{Step 2: second homotopy}
\\[1ex]
We consider the function $F:\ V\times [0,1] \rightarrow \xR^N\times \xR^M\times \xR^M$ given by:
\begin{equation}
\begin{array}{ll}
F(u,p,z,\theta)=
\left| \begin{array}{l} \displaystyle
v_{\edge,i}=
a(u,\varphi_\edge^{(i)})
- \theta \int_{\Omega,h} p\ \nabla \cdot \varphi_\edge^{(i)} \, {\rm d}x- \int_\Omega f \cdot \varphi_\edge^{(i)} \, {\rm d}x,
\hfill
\edge \in \edgesint,\ 1 \leq i \leq d
\\ [2ex] \displaystyle
q_K=\frac{|K|}{\dt} (p_K-p_K^\ast) 
+ \theta \sum_{\edge=K|L} \vsp\, p_K - \vsm\, p_L,
\hspace{20ex}
K \in \mesh
\\ [2ex] \displaystyle
s_K=\frac{|K|}{\dt} (z_K- \frac{p_K^\ast}{\rho_K^\ast}\, z_K^\ast) + \theta \sum_{\edge=K|L} \vsp\, z_K - \vsm\, z_L,
\hfill
K \in \mesh
\end{array} \right .
\end{array}
\label{eqalpha2} \end{equation}
The system $F(u,p,z,1)=0$ is the same as the system obtained at the end of the preceding homotopy for $\theta=0$, and the system $F(u,p,z,0)=0$ is linear and clearly regular (by stability of the bilinear form $a(\cdot,\cdot)$).

\medskip
In addition, the third equation is now decoupled from the first two ones, and these latter have the structure of a monophasic compressible problem as studied in \cite{gal-07-an}.
From this theory, an estimate similar to the first one in the preceding step is available and reads:
\[
\frac 1 2 \normad{u} + \frac 1 \dt \sum_{K\in\mesh} |K|\ p_K\ \log(p_K) \leq 
\frac 1 \dt \sum_{K\in\mesh} |K|\ p^\ast_K\ \log(p^\ast_K) + \frac 1 2 \normmad{f}
\]
Since the function $s \mapsto s \log(s)$ is bounded by below on $(0,+\infty)$, this latter relation yields:
\begin{equation}
\norma{u} \leq c_{10}
\label{est_u2}\end{equation}
By lemma \ref{bound_rho_z}, we thus directly get:
\begin{equation}
c_{11} \leq p \leq c_{12},\qquad c_{13} \leq z \leq c_{14}
\label{est_uz2}\end{equation}

\vspace{3ex}
\noindent \underline{Conclusion}
\\[1ex]
We choose $\epsilon$ small enough for the relations \eqref{est_u}, \eqref{est_z_inf}, \eqref{est_p_sup_prel} and \eqref{est_p_sup} to hold,
$\epsilon < \min(c_{11},\, c_{13})$ and, in addition:
\[
\epsilon < \max(c_3,\, c_5) \ |\log(\epsilon)|^{-1/2}
\]
which is possible because the function $s \mapsto s \log s$ tends to zero when $s$ tends to zero.
Let now $M$ be such that:
\[
M > \max \left[ \max(c_1,c_9) \ |\log (\epsilon)|^{1/2},\ c_4,\, c_{10},\, c_{12},\, c_{14} \right]
\]
Then, from inequalities \eqref{est_u}, \eqref{est_z_inf}, \eqref{est_z_sup}, \eqref{est_p_inf}, \eqref{est_p_sup}, \eqref{est_u2} and \eqref{est_uz2}, we get that throughout both homotopies, the unknown $(u,p,z)$ remains in ${\cal O}$.
As the last linear system is regular and admits a solution in ${\cal O}$, the topological degree of $F(\cdot,\cdot,\cdot,\theta)$ with respect to ${\cal O}$ and zero remains different of zero all along both homotopies, which proves the existence of a solution in ${\cal O}$.

\medskip
We now turn to the proof of the {\em a priori} estimates $\rho>0$, $z>0$, $0< y^\ast \leq 1$ and $p>0$.
The fact that, if $\rho^\ast>0$ and $z^\ast>0$, then $\rho>0$ and $z>0$ is a direct consequence of lemma \ref{bound_rho_z} applied to the second and third relation of problem \eqref{pbdiph_disc}.
In addition, as both $\rho^\ast>0$ and $\rho>0$, lemma \ref{larrout'} applies and thus, as $0< y^\ast \leq 1$, we have $0< y \leq 1$.
If $p \geq \rho_\ell$, the fact that $p>0$ is evident.
In the other case, by the equation of state written as a function of $p$ and $y$ (third form of \eqref{rhomixture}), we get first that:
\begin{equation}
p \leq \rho < \rho_\ell
\label{rho=rhomoy}\end{equation}
and, second, that, since $\rho >0$, the pressure does not vanish.
The second form of this same relation \eqref{rhomixture} thus can be written:
\[
\rho= \frac z p \, p + (1- \frac z p)\, \rho_\ell=\alpha_g \, p + (1-\alpha_g) \, \rho_\ell
\]
As $p\neq \rho_\ell$, the void fraction $\alpha_g$ thus reads:
\[
\alpha_g=\frac{\rho-\rho_\ell}{\rho_\ell - p}
\]
which, by inequalities \eqref{rho=rhomoy}, yields $\alpha_g>0$ and, finally, since $z>0$, $p>0$.
\end{proof}

This existence result applies directly to the pressure correction step used in the algorithm presented in this paper, with a particular expression for the bilinear form $a(\cdot,\cdot)$, which reads, dropping for short the time exponents:
\[
a(u,v)=\sum_{\edge \in \edgesint} \frac{|D_\edge|} \dt\ \rho_\edge \ u_\edge \cdot v_\edge
\]
Note that the analysis is performed here with a very simple equation of state for the gas ($p=\rho$), but would be readily extended to general barotropic laws $p=\wp(\rho)$, under the mild assumptions that the corresponding free energy exists and is convex and the function $\wp(\cdot)$ is increasing and one to one from $(0,+\infty)$ to $(0,+\infty)$.

\medskip
Let us now turn to the discretization of a stationary diphasic problem.
As happens in the monophasic case, \cite{gal-07-conv}, it is likely that, in the case where the velocity is prescribed on the whole boundary, this problem needs to be completely determined the data of the total mixture mass (say $M_{m}$) and of the total gas mass (say $M_g$) present in the computational domain.
A natural way to impose these two conditions is to add to this problem two regularizing terms in the mass balance and the gas mass balance:
\[
\left| \begin{array}{ll} \displaystyle
c(h)\ |K|\ \left[\varrho^{\,p,z}(p_K,\, z_K)-\frac{M_m}{|\Omega|}\right] + \sum_{\edge=K|L} \vsp\ \varrho^{\,p,z}(p_K,z_K) - \vsm\ \varrho^{\,p,z}(p_L,z_L)=0
\quad & \displaystyle
\forall K \in \mesh
\\[2ex] \displaystyle
c(h)\ |K|\ \left[z_K-\frac{M_g}{|\Omega|}\right] + \sum_{\edge=K|L} \vsp\ z_K - \vsm\ z_L=0
\quad & \displaystyle
\forall K \in \mesh
\end{array} \right.
\]
where $c(h)$ is a regularization parameter tending to zero with the size of the mesh.
In this case, the present existence theory directly applies, provided that the momentum balance equation remains linear with respect to the velocity.
Of course, under the same restriction, this is true also for an implicit discretization of a time-dependent problem.

\medskip
In view of the stability results provided for the advection operator, adding such a term to the first relation of the problem (\ie\ the momentum balance) should lead to a rather straightforward extension of the present existence result; the advection term would be multiplied by the homotopy parameter and the stability (\ie\ an analogue to estimate \eqref{est_u}) would stem from the diffusion term.
Note that, in this case, to keep the stability of the advection term, a regularization term consistent with the mass balance one should also be introduced in the momentum balance equation.

\medskip
We have shown in this paper that, with a Darcy's law for the drift velocity and a particular discretization for this term, the drift term is dissipative.
So this term does not prevent to obtain stability estimates as \eqref{est_u}; this suggests that the existence theory developped here may perhaps be extended to the complete drift flux model.

\medskip
Finally, we have not dealt in this study with the case where liquid monophasic zones ($z=0$) exists in the flow.
In such zones, the pressure changes of mathematical nature: it is no more a parameter entering the equation of state and determined by the local density, but a Lagrange multiplier for the incompressibility constraint.
Note that this fact is already underlying in the present study: indeed, the incompressibility of the liquid prevents to derive $\xLinfty$ estimates for the pressure from $\xLinfty$ estimates for the density (which are readily obtained using a conservation argument), and we must invoke to this purpose the stability of the discrete gradient (\ie\ the discrete {\em inf-sup} condition), that is typically the argument allowing to control the pressure in incompressible flow problems.
However, obtaining {\it a priori} estimates when $z$ may vanish in the flow seems a difficult task, which should deserve more efforts.
On the contrary, obtaining existence results for two barotropic phases seems to be rather simpler than the analysis performed here.


\begin{acknowledgement}{\em Acknowledgements.} The authors thank F. Duval, from IRSN, for helpful discussions in the course of this work, and F. Babik, from the ISIS development team at IRSN, for supporting the implementation of this scheme.
\end{acknowledgement}

\bibliographystyle{plain}
\bibliography{./mainbib}
\end{document}